\documentclass[twoside,12pt,leqno]{article}
\advance\topmargin by -\headheight\advance\topmargin by -\headsep
\newcommand{\tmvolxx}{xx}
\newcommand{\tmyearyyyy}{yyyy}

\newcommand{\FirstPageHead}[3]{
{\footnotesize 
\vskip -8mm 
\centerline {Travaux math\'ematiques, \quad 
Volume #1 (#2), 
#3,\quad \copyright\  Universit\'e du Luxembourg}}\vspace{-3mm}}

\usepackage{amsmath}
\usepackage{amsthm}
\usepackage{amssymb}
\usepackage{amscd}
\usepackage{graphicx}
\usepackage{epsfig}
\numberwithin{equation}{section}
\newtheorem{theorem}{Theorem}[section]
\newtheorem{lemma}[theorem]{Lemma}
\newtheorem{proposition}[theorem]{Proposition}

\newtheorem{conjecture}[theorem]{Conjecture}
\theoremstyle{definition}
\newtheorem{definition}[theorem]{Definition}

\newtheorem{remark}[theorem]{Remark}

\numberwithin{equation}{section}

\allowdisplaybreaks

\usepackage{slashed}
\usepackage[toc,page]{appendix}
\usepackage{tikz}

\begin{document}
\thispagestyle{empty}
\FirstPageHead{\tmvolxx}{\tmyearyyyy}{\pageref{firstpage}--\pageref{lastpage}}
\label{firstpage}
\newcommand{\df}{\displaystyle\frac}    
\newcommand{\seq}[1]{\left<#1\right>}  

\newcommand{\boldnabla}{\mbox{\boldmath$\nabla$}}
\newcommand{\tildeboldnabladelta}{\tilde \boldnabla}
\newcommand{\wf}{\operatorname{WF}}
\newcommand{\Span}{\operatorname{span}}
\newcommand{\Hom}{\operatorname{H}}
\newcommand{\PSL}{\operatorname{PSL}}
\newcommand{\PSLC}{\operatorname{PSL}(2\text{, }\Complex)}
\newcommand{\Homo}{\mathop{\fam0 Hom}\nolimits}
\newcommand{\Diffo}{\operatorname{Diff}^+}
\newcommand{\SqInt}{\operatorname{L}^2}
\newcommand{\End}{\mathop{\fam0 End}\nolimits}
\newcommand{\Alt}{\text{Alt} }                                       
\newcommand{\Ker}{\text{Ker} }  
\newcommand{\im}{\mathop{\fam0 Im}\nolimits}
\newcommand{\re}{\mathop{\fam0 Re}\nolimits}
\newcommand{\diff}{\mathrm{d}}
\newcommand{\LN}{\operatorname{L}_N}
\newcommand{\CP}{\mathop{\fam0 {\mathbb CP}}\nolimits}
\newcommand{\Half}{\mathbb H}
\newcommand{\Hyp}[1]{\mathbb{H}^#1}
\newcommand{\apb}{\tilde{a}_+}
\newcommand{\bpb}{\tilde{b}_+}
\newcommand{\cpb}{\tilde{c}_+}
\newcommand{\sapb}{\tilde{a}_+\sqrt{N}}
\newcommand{\sbpb}{\tilde{b}_+\sqrt{N}}
\newcommand{\scpb}{\tilde{c}_+\sqrt{N}}
\newcommand{\amb}{\tilde{a}_-}
\newcommand{\bmb}{\tilde{b}_-}
\newcommand{\cmb}{\tilde{c}_-}
\newcommand{\samb}{\tilde{a}_-\sqrt{N}}
\newcommand{\sbmb}{\tilde{b}_-\sqrt{N}}
\newcommand{\scmb}{\tilde{c}_-\sqrt{N}}
\newcommand{\rhoq}{\eta}
\newcommand{\psitp}{\tilde{\psi}^{\prime}}
\newcommand{\cth}{c_{\bpar}}
\newcommand{\Dth}{D_{\bpar}}
\newcommand{\conj}[1]{\overline{#1}}
\newcommand{\qPoch}[2]{\left(#1 \text{;} #2\right)_{\infty}}
\newcommand{\Real}{\mathbb{R}}
\newcommand{\MCG}{\operatorname{MCG}_{g,s}}
\newcommand{\Teich}{\mathcal{T}_{g,s}}
\newcommand{\DTeich}{\widetilde{\mathcal{T}}_{g,s}}
\newcommand{\Ratio}{{\mathcal{R}}(\tau)}
\newcommand{\RatioC}{{\mathcal{R}}_{\Complex}(\tau)}
\newcommand{\RatioWell}{\mathcal{R}(\GenSurf)}
\newcommand{\RatioCWell}{\mathcal{R}_{\Complex}(\GenSurf)}
\newcommand{\GenSurf}{\Sigma_{g\text{,}s}}
\newcommand{\hol}{\mathop{\fam0 Hol}\nolimits}
\newcommand{\PtGr}{\mathcal{G}(\GenSurf)}
\newcommand{\glie}{\mathfrak g}
\newcommand{\glc}{\mathfrak g_{\Complex}}
\newcommand{\Shp}{\operatorname{S}}
\newcommand{\GShp}{\widetilde{\Shp}}
\newcommand{\LShp}{\operatorname{LS}}
\newcommand{\LGShp}{\widetilde{\LShp}}
\newcommand{\dfield}[1]{\frac{\partial}{\partial #1}}
\newcommand{\D}{\mathcal{D}}
\newcommand{\abs}[1]{\left\lvert #1 \right\rvert}
\newcommand{\Ha}{\operatorname{H}_a}
\newcommand{\Integer}{\mathbb{Z}}
\newcommand{\IntN}{\Integer /N\Integer}
\newcommand{\Complex}{\mathbb{C}}
\newcommand{\Affine}{\mathbb{A}}
\newcommand{\Torus}{\mathbb{T}}
\newcommand{\Bordism}{\mathcal{B}}
\newcommand{\Hilb}{ {\mathcal H}}
\newcommand{\Schw}{\mathcal{S}}
\newcommand{\Line}{\mathcal{L}}
\newcommand{\Cont}{\operatorname{C}}
\newcommand{\Smooth}{\Cont^{\infty}}
\newcommand{\OBJ}{\mathrm{Obj}}
\newcommand{\morph}{\mathrm{Mor}}

\newcommand{\sign}{\operatorname{sign}}
\newcommand{\bpar}{\mathrm{b}}
\newcommand{\cb}{c_{\bpar}}
\newcommand{\bdilog}{\Phi_{\bpar}}
\newcommand{\bDilog}{\mathrm{D}_{\bpar}}
\newcommand{\cthe}{c_{\bpar}}
\newcommand{\thdilog}{\Phi_{\bpar}}
\newcommand{\thDilog}{\mathrm{D}_{\bpar}}
\newcommand{\SthDilog}{\slashed{\mathrm{D}}_{\bpar}}
\newcommand{\thPsylog}{\psi_{\bpar}}
\newcommand{\Psylog}{\mathsf{\Psi}_{\bpar}}
\newcommand{\opPsylog}{\mathsf{\Psi}}
\newcommand{\SthPsylog}{\slashed{\Psi}_{\bpar}}
\newcommand{\Poch}[3]{\left(#1 \text{;} #2\right)_{#3}}
\newcommand{\Dilog}{\operatorname{Li}_2}

\newcommand{\bkt}[1]{\langle #1\rangle}
\newcommand{\bra}[1]{\langle #1 |}
\newcommand{\ket}[1]{| #1\rangle}
\newcommand{\Reala}[1]{\langle #1 |}
\newcommand{\comm}[1]{\left[ #1\right]}
\newcommand{\opq}{\mathsf{q}}
\newcommand{\opp}{\mathsf{p}}
\newcommand{\opu}{\mathsf{u}}
\newcommand{\opv}{\mathsf{v}}
\newcommand{\opw}{\mathsf{w}}
\newcommand{\opT}{\mathsf{T}}
\newcommand{\opG}{\mathsf G}
\newcommand{\opS}{\mathsf{S}}
\newcommand{\opF}{\mathsf{F}}
\newcommand{\opL}{\mathsf{L}}
\newcommand{\opD}{\mathsf{D}}
\newcommand{\opx}{\mathsf{x}}
\newcommand{\opy}{\mathsf y}
\newcommand{\opl}{\mathsf l}
\newcommand{\opm}{\mathsf m}
\newcommand{\opr}{\mathsf r}
\newcommand{\opb}{\mathsf b}
\newcommand{\opj}{\mathsf j}
\newcommand{\opA}{\mathsf A}
\newcommand{\opB}{\mathsf B}
\newcommand{\opf}{\mathsf f}
\newcommand{\opg}{\mathsf g}
\newcommand{\p}{{\prime}}
\newcommand{\pp}{{\prime\prime}}




\markboth{J.E. Andersen, S. Marzioni}{Teichm\"{u}ller TQFT at level $N$}
$ $
\bigskip

\bigskip

\centerline{{\Large   Level $N$ Teichm\"{u}ller TQFT and}} \vskip 0.5cm

\centerline{{\Large Complex Chern--Simons Theory}}

\bigskip
\bigskip
\centerline{{\large by  J{\o}rgen Ellegaard Andersen and Simone Marzioni\footnote{Work supported in part by the center of excellence grant ``Center for Quantum Geometry of Moduli Spaces'' from the Danish National Research Foundation (DNRF95).}
}}

\vspace*{.7cm}

\begin{abstract}
	In this manuscript we review the construction of the Teichm\"{u}ller TQFT in \cite{AK1}, upgrading it to a theory dependent on an extra odd integer $N$ using results developed in \cite{AK3}. We also describe how this theory is related with quantum Chern--Simons Theory at level $N$ with gauge group $\PSLC$. 
\end{abstract}

\pagestyle{myheadings}
\section{Introduction}

In this paper we review Andersen and Kashaev's construction of the Teichm{\"u}ller TQFT from \cite{AK1} making it dependent on an extra odd integer $N$, called the level. The original work of \cite{AK1} corresponds to the choice $N=1$, and emerged as an extension to a $3$--dimensional theory of the representations one obtains from Quantum Teichm\"uller Theory  \cite{KQuantTeich}. In particular, it defines a class of quantum invariants for hyperbolic knots, dependent on a continuous parameter $\bpar$.   The level $N$ Teichm\"uller TQFT is an analogously upgarde of representations in  Quantum Teichm\"uller theory and it depends on a pair of parameters $\left(\bpar, \, N\right)$, one continuous and one discrete. Such quantum theory is related to the level $N$ Chern--Simons theory with gauge group $\PSLC$ via the level $N$ Weil-Gel'fand-Zak transform. Such a relation was proposed in \cite{AK3}, and here we show it in a more tight way for the four punctured sphere.
One of the main ingredient in the construct of the Teichm\"uller TQFT is the \emph{quantum dilogarithm}, that is a function $\thDilog \colon\Real\times\IntN\rightarrow \Complex$, satisfying some particular properties. Such functions were introduced in \cite{AK3}, which for $N=1$ is Faddeev's original quantum dilogarithm. The theory we get has different and interesting unitarity behaviour depending on the pair of parameters $(\bpar, N)$: for level $N=1$ the theory is unitary  whenever $\bpar>0$ or $\abs{\bpar} =1$ while for higher level $N>1$ the unitarity is only manifest when $\abs{\bpar}=1$ while in the case $\bpar>1$ the behaviour is more exotic. We will consider both situations here and we will use the setting $\bpar >0$ to present some asymptotic properties in the limit $\bpar \rightarrow 0$.
The Teichm\"uller TQFT can be used to define knot invariants starting from triangulations of their complements. In this presentation we update the examples presented in \cite{AK1} to the level $N$ setting together with their asymptotic analysis. For the simplest hyperbolic knot we show the appearance of the Baseilhac--Benedetti invariant from \cite{BBqhi} in such a limit.

It is an interesting challenge how the TQFT's which are reviewed in this paper are related to the Witten-Reshetikhin-Turaev TQFT's \cite{W1, RT1,RT2, BHMV1,BHMV2, B1} and in particular how they are related to the geometric construction of these TQFT's \cite{ADW,H1,La1,TUY,AU1,AU2,AU3,AU4} and to Witten's proposal for the construction of the complex quantum Chern-Simons theory \cite{W3}, which can actually be constructed from a purely mathematical point of view \cite{AG}, resulting also in the mathematically well-defined Hitchin-Witten connection in the bundle of quantizations of the moduli space of flat $\text{SL}(n,{\mathbb C})$-connections
over Teichm\"{u}ller space. In the classical case of compact groups, the description of the representations of the mapping class groups via the monodromy of the Hitchin connections turned out to be very useful to prove deep properties about these representations \cite{A1,A2, A4, Aqt,AH,  AHJMMc}, some of which also uses the theory of Toeplitz operators \cite{BMS,KS}. Understanding how these kinds of results can be extended to the complex quantum theory discused in this paper will be very interesting and most likely involve using Higgs bundles techniques \cite{H2}. Certainly we have already seen the start of this with the Verlinde formula for Higgs bundle moduli space \cite{AGP}.

The paper is organised as follows.
In section \ref{sbs:PtRep} we recall the definition of the (decorated) Ptolemy groupoid of punctured surfaces, which is the combinatorial foundation over which Quantum Teichm\"uller Theory is defined.
In  section \ref{sc:qDilog} we recall the quantum dilogarithm $\thDilog$, and we list some of its properties. The function $\thDilog$ was introduced in \cite{AK3} for the first time, but some of its properties that we list here are not present in the literature.
In section \ref{sc:ModelSpace} we carry out the quantization of the moduli space of $\PSLC$ flat connections over a four punctured sphere with unipotent holonomies around the punctures. We follow the prescriptions of geometric quantisation,  together with a choice of real polarisation, and we connect the resulting algebra of observables to the $\SqInt(\Real\times\IntN)$ representations of quantum Teichm\"uller theory of the previous section.
Finally in section \ref{ch:TQFT} we construct the the Teichm\"uller TQFT functor $F_\bpar^{(N)}$ mirroring the construction in \cite{AK1} and we study some examples and their asymptotic behaviour. It would of course be interesting to go through all the examples treated in \cite{AN} in this volume for the level $N$ theory as well.

\subsection*{Acknowledgements} We would like to thank Rinat Kashaev for many interesting discussions.

\section{Ptolemy Groupoid}\label{sbs:PtRep}

Let $\GenSurf$ be a surface of genus $g$ with $s$ punctures, such that $s>0$ and $2 - 2g + s <0$.
\begin{definition}
	An \emph{ideal arc} $\alpha$ is the homotopy class relative endpoints of the embedding of a path in $\Sigma_{g,s}$, such that the endpoints are punctures of the surface.\\
	An \emph{ideal triangle} is a triangle with the vertices removed, such that the edges are ideal arcs.\\
	An \emph{ideal triangulation} $\tau$ of $\Sigma_{g,s}$ is a collection of disjoint ideal arcs such that $\Sigma_{g,s} \setminus \tau$ is a collection of interiors of ideal triangles.
\end{definition}
Given an ideal triangulation $\tau$, $\Delta_j(\tau)$ will denote the set of its $j$-dimensional cells.
\begin{definition}
	A \emph{decorated ideal triangulation}  of $\GenSurf$ is an ideal triangulation $\tau$ up to isotopy relative to the punctures, together with the choice of a distinguished corner in each ideal triangle and a bijective ordering map 
	$$\overline{\tau} : \{1\text{, }\dots\text{, } s\} \ni j \mapsto \overline{\tau}_j \in \Delta_2(\tau)\text{.}$$
	We denote the set of all decorated ideal triangulation as $\dot{\Delta} = \dot{\Delta}(\GenSurf)$.
\end{definition}
When we say that $\tau$ is a decorated ideal triangulation we mean that $\tau$ is the set of decorated ideal triangles in the triangulation. 

One of the main interests in quantizing moduli spaces is the consequent construction of representations of (central extensions of) the mapping class group of the surfaces \cite{W3, H2, La1, A1, AU4, AG}. Quantum Teichm\"uller theory produce instead representations of a bigger object  called the (decorated) \emph{Ptolemy Groupoid} that we are going to introduce now.

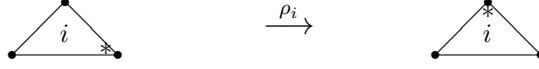
\begin{figure}[htb]
	\centering
	\begin{picture}(200,20)
	\put(0,0){\begin{picture}(40,20)
		\put(0,0){\line(1,0){40}}
		\put(0,0){\line(1,1){20}}
		\put(20,20){\line(1,-1){20}}
		\put(0,0){\circle*{3}}
		\put(20,20){\circle*{3}}
		\put(40,0){\circle*{3}}
		\footnotesize
		\put(33,0){$*$}
		\put(18,5){$i$}
		\end{picture}}
	\put(160,0){\begin{picture}(40,20)
		\put(0,0){\line(1,0){40}}
		\put(0,0){\line(1,1){20}}
		\put(20,20){\line(1,-1){20}}
		\put(0,0){\circle*{3}}
		\put(20,20){\circle*{3}}
		\put(40,0){\circle*{3}}
		\footnotesize
		\put(17.5,14){$*$}
		\put(18,5){$i$}
		\end{picture}}
	\put(95,8){$\stackrel{\rho_i}{\longrightarrow}$}
\end{picture}
\caption{Transformation $\rho_i$.}\label{fg:Pt1}
\end{figure}

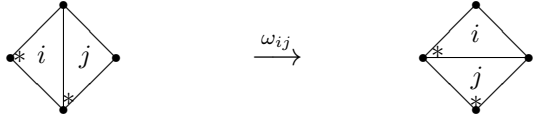
\begin{figure}[htb]
	\centering
	\begin{picture}(200,40)
	\put(0,0){
		\begin{picture}(40,40)
		\put(20,0){\line(-1,1){20}}
		\put(40,20){\line(-1,-1){20}}
		\put(0,20){\line(1,1){20}}
		\put(40,20){\line(-1,1){20}}
		\put(20,0){\line(0,1){40}}
		\put(20,0){\circle*{3}}
		\put(0,20){\circle*{3}}
		\put(20,40){\circle*{3}}
		\put(40,20){\circle*{3}}
		\footnotesize
		\put(10,18){$i$}\put(26,18){$j$}
		\put(1,18){$*$}
		\put(19.5,2){$*$}
		\end{picture}}
	\put(160,0){\begin{picture}(40,40)
		\put(20,0){\line(-1,1){20}}
		\put(40,20){\line(-1,-1){20}}
		\put(0,20){\line(1,1){20}}
		\put(40,20){\line(-1,1){20}}
		\put(0,20){\line(1,0){40}}
		\put(20,0){\circle*{3}}
		\put(0,20){\circle*{3}}
		\put(20,40){\circle*{3}}
		\put(40,20){\circle*{3}}
		\footnotesize
		\put(18,26){$i$}\put(18,10){$j$}
		\put(3,20){$*$}
		\put(17.5,1){$*$}
		\end{picture}}
	\put(95,17){$\stackrel{\omega_{ij}}{\longrightarrow}$}
\end{picture}
\caption{Transformation $\omega_{ij}$.}\label{fg:Pt2}
\end{figure}

Recall that, given a group $G$ acting freely on a set $X$, we can define an associated groupoid $\mathcal{G}$ as follows. The objects of $\mathcal{G}$ are $G$-orbits in $X$ while morphisms are $G$-orbits in $X\times X$ with respect to the diagonal action. Then for any $x\in X$ we can consider the object $[x]$ and for any pair $(x,y)\in X\times X$ we can consider the morphism  $[x,y]$. When $[y] = [u]$ there will be a $g\in G$ so that $gu = y$ and we can define the composition $[x,y][u,v] = [x,gv]$. The unit for $[x]$ is given by $[x,x]$. If the action of $G$ is transitive, we would get an actual group.
We will abbreviate $[x_1,x_2][x_2,x_3] \cdots [x_{n-1},x_n]$ with $[x_1,x_2,\dots , x_n]$.\\
We define the \emph{decorated Ptolemy groupoid} $\PtGr$ of a punctured surface $\GenSurf$ following the above recipe. The set we consider is the set  $\dot{\Delta}$ of decorated triangulations $\tau$ of $\GenSurf$. The  free group action is the one of the mapping class group $\MCG$ acting on $\dot \Delta$. This action is not transitive, meaning that not all pairs of decorated ideal triangulations can be related by a mapping class group element. However in the language of groupoids, we can still describe generators and relations for the morphism groups. 
For $\tau\in\dot{\Delta}$ there are three kind of generators $[\tau,\tau^{\sigma}]$, $[\tau, \rho_i\tau]$ and $[\tau, \omega_{i,j}\tau]$, where $\tau^{\sigma}$ is obtained by applying the permutation $\sigma\in\mathbb{S}_{|\tau |}$ to the ordering of triangles in $\tau$, $\rho_i\tau$ is obtained by changing the distinguished corner in the triangle $\bar{\tau}_i\in\tau$ as in Figure \ref{fg:Pt1}, and $\omega_{i,j}$ is obtained by applying a decorated diagonal exchange to the quadrilateral  made of the two decorated ideal triangles $\bar{\tau}_i$  and $\bar{\tau}_j$ as in Figure \ref{fg:Pt2}.

The relations are usually grouped in two sets, the first being
\begin{equation}
[\tau,\tau^{\alpha}, (\tau^{\alpha})^{\beta}] = [\tau,\tau^{\alpha\beta}] \text{,   }\quad \alpha,\beta\in\mathbb{S}_{\tau}\text{,}
\end{equation}
\begin{equation}
[\tau,\rho_i\tau,\rho_i\rho_i\tau, \rho_i\rho_i\rho_i\tau] =\mathrm{id}_{[\tau]} \text{,}
\end{equation}
\begin{equation}\label{eq:Pent1}
[\tau,\omega_{i,j}\tau,\omega_{i,k}\omega_{i,j}\tau,\omega_{j,k}\omega_{i,k}\omega_{i,j}\tau] =[\tau,\omega_{j,k}\tau,\omega_{i,j}\omega_{j,k}\tau] 
\end{equation}
\begin{equation}\label{eq:Inv1}
[\tau,\omega_{i,j}\tau,\rho_i\omega_{i,j}\tau,\omega_{j,i}\rho_i\omega_{i,j}\tau] =[\tau,\tau^{(i,j)}, \rho_j\tau^{(i,j)}, \rho_i\rho_j\tau^{(i,j)}] 
\end{equation}
The first two relations are obvious, the third is called the Pentagon Relation and  the fourth is  called the Inversion Relation.\\

The second set of relations, are commutation relations
\begin{equation}
[\tau,\rho_i\tau, \rho_i\tau^{\sigma}] = [\tau,\tau^\sigma, \rho_{\sigma^{-1}(i)}\tau^\sigma] \text{,}
\end{equation}
\begin{equation}
[\tau,\omega_{i,j}\tau, (\omega_{i,j}\tau)^{\sigma}] = [\tau,\tau^\sigma, \omega_{\sigma^{-1}(i)\sigma^{-1}(j)}\tau^\sigma] \text{,}
\end{equation}
\begin{equation}
[\tau,\rho_j\tau, \rho_j\rho_i\tau] = [\tau,\rho_i\tau, \rho_i\rho_j\tau] \text{,}
\end{equation}
\begin{equation}
[\tau,\rho_i\tau, \omega_{j,k}\rho_i\tau] = [\tau,\omega_{j,k}\tau, \rho_i\omega_{j,k}\tau] \text{,   } i\notin \{j,k\}\text{,} 
\end{equation}
\begin{equation}
[\tau,\omega_{i,j}\tau, \omega_{k,l}\omega_{i,j}\tau] = [\tau,\omega_{k,l}\tau, \omega_{i,j}\omega_{k,l}\tau] \text{,   } \{ i,j\}\cap\{ k,l\} =\emptyset \text{,} 
\end{equation}

To every decorated ideal triangulation $\tau\in\dot{\Delta}$ it is possible to associate a simple symplectic space $\Ratio$, called the space of \emph{Ratio Coordinates}.
We summarize here its relation with the Ptolemy groupoid and refer to \cite{KQuantTeich} for a detailed introduction to ratio coordinates and their realtion to the Teichm\"uller space. 
Let $M \equiv 2g - 2 + s = \abs{\tau}$ be the number of ideal triangles, then $\Ratio \equiv \left(\Real_{>0}\times\Real_{>0}\right)^M$. 
Let $x^j\equiv(x^j_1,x^j_2)\in \Real_{>0}\times\Real_{>0}$, for $j= 1,\,\dots,M$ be the coordinates associated to the ideal triangle $\bar{\tau}_j\in\Delta_2(\tau)$. 
The symplectic form that we consider on $\Ratio$ is 
\begin{align}
\omega_\tau \equiv \sum_{j=1}^M \frac{\diff x^j_1}{x^j_1}\wedge \frac{\diff x^j_2}{x^j_2}
\end{align}
Now we want to describe the action of $\PtGr$ as symplectomorphisms between  these spaces.  The morphisms $[\tau, \tau^{\sigma}]$ act by permuting the coordinates in $\Ratio$. The  morphism $[\tau, \rho_i\tau]$ acts as the identity on any pair $x=(x_1,x_2)$ corresponding to ideal triangles different from $\bar{\tau}_i$ and as $(x_1,x_2)=x \mapsto y = (\frac{x_2}{x_1} , \frac{1}{x_1})$ for the pair of coordinates corresponding to $\bar{\tau_i}$. Finally the action of $[\tau, \omega_{i,j}\tau]$ is the identity on $\bar{\tau}_k$ for $k \neq i$, $j$ while letting $x=(x_1,x_2)$ and $y=(y_1,y_2)$ be the coordinates corresponding to the triangles $\bar{\tau}_i$ and $\bar{\tau}_j$ respectively, and letting  $u=(u_1,u_2)$ and $v= (v_1, v_2)$ be the coordinates of the triangles $\overline{\omega_{i,j}\tau}_i$ and $\overline{\omega_{i,j}\tau}_j$, then we have $u= x\bullet y$ and $v= x * y$ where 
\begin{align}\label{eq:RatioPtolemy}
&x\bullet y := \left(x_1y_1, x_1y_2 + x_2\right)\\
\nonumber
&x*y := \left(\frac{y_1x_2}{x_1y_2 +x_2}, \frac{y_2}{x_1y_2 + x_2}\right)\text{.}
\end{align}

Let $\widetilde{\Delta}$ be the set of pairs $(\tau,\Ratio)$, $\tau\in\dot\Delta$. Then the space $\RatioWell$  is defined as the quotient of $\widetilde{\Delta}$  by the action of $\PtGr$ as described above. This space of coordinates is now independent of the triangulation. For more details on the (decorated or not) Ptolemy groupoid see \cite{Pbook},\cite{FKmcg},\cite{KLiuTeich2}\cite{KDehnDil1}.

\section{Quantum Dilogarithm}\label{sc:qDilog}

In this section we recall the  quantum dilogarithm $\thDilog$ over $\Affine_N$ and we state some of their properties.  
\begin{definition}[$q$-Pochammer Symbol]\label{df:Poch}
	Let $x, q\in\Complex$, such that  $|q| <1$. Define the \emph{$q$-Pochammer Symbol} of $x$ as 
	$$\qPoch{x}{q} := \prod_{i=0}^{\infty} (1-xq^i)$$
\end{definition}

\begin{theorem}\label{th:PochPent}
	Let $X$, $Y$  satisfying $XY= qYX$. Then the following five-term relation holds true
	\begin{equation}
	\qPoch{Y}{q}\qPoch{X}{q} = \qPoch{X}{q}\qPoch{-YX}{q}\qPoch{Y}{q}\text{.}
	\end{equation} 
\end{theorem}

\begin{definition}[Faddeev's Quantum Dilogarithm \cite{Fqdilog}] \label{df:Faddeev}
	Let $z, \bpar\in\Complex$ be such that $\re\bpar\neq 0$,  $|\im(z)|< |\im(\cb)|$, where $\cb := i(\bpar  + \bpar^{-1})/2$. Let $C\subset\Complex$, $C= \Real +i0$ be a contour equal to the  the real line outside a neighborhood  of the origin that avoid the singularity in $0$ going in the upper half plane.  \emph{Faddeev's quantum dilogarithm} is defined to be 
	\begin{equation}
	\bdilog (z) = \text{exp} \left(\int_{C} \frac{e^{-2izw}\diff w}{4\sinh (w\bpar) \sinh(w\bpar^{-1})w} \right)\text{.} 
	\end{equation}
\end{definition}
It is evident that $\bdilog$ is invariant under the following changes of parameter
\begin{equation}
\bpar\, \leftrightarrow\, \bpar^{-1}\,\leftrightarrow\, -\bpar\text{,}
\end{equation}
so that our choice of $\bpar$ can be restricted to the first quadrant
\begin{equation}
\re\bpar >0\text{, }\qquad \im\bpar \geq 0
\end{equation}
which implies 
\begin{equation}
\im(\bpar^2) \geq 0.
\end{equation}

Faddeev's quantum dilogarithm has a lot of other interesting properties and applications, see for example \cite{Fqdilog},\cite{FKqdilog},\cite{FKV1} and \cite{Vhypergeometry}.

Let $N \geq 1$ be a positive \emph{odd} integer. 
Then, following \cite{AK3}, we can define a \emph{quantum dilogarithm} over $\Affine_N$ as follows

\begin{equation}\label{df:thDilog}
\thDilog (x,n) := \prod_{j=0}^{N-1} \thdilog\left(\frac{x}{\sqrt{N}} + (1-N^{-1})\cthe -i\bpar^{-1}\frac{j}{N} -i\bpar \left\{ \frac{j+n}{N}\right\}\right)
\end{equation}
where $\{p\}$ is the fractional part of $p$, and $\thdilog$ is the Faddeev's quantum dilogarithm. Of course for $N=1$ we have just $\thdilog(x)$.
The function $\thDilog$ was introduced in \cite{AK3} only for $\abs{\bpar} = 1$. It satisfies a series properties that we are going to list.
\begin{lemma}[Inversion Relation \cite{AK3}]\label{lm:thDilogInv}
	$$\thDilog(x,n)\thDilog(-x,-n) = e^{\pi i x^2} e^{-\pi in(n+N)/N}\zeta_{N,\, inv}^{-1} \text{,}$$
	where
	\begin{align}
	\zeta_{N,\, inv} = e^{\pi i (N+2\cthe^2 N^{-1})/6}\text{.}
	\end{align}
\end{lemma}
Unitarity properties are different in the two situations $\abs{\bpar} = 1$ or $\bpar\in\Real$.
\begin{lemma}[Unitarity]\label{lm:thDilogUnit}
	\begin{align}
	\overline{\thDilog(x,n)} = \thDilog(\bar x,n)^{-1} \qquad &\text{ if } \abs{\bpar} = 1\text{,}\\
	\overline{\thDilog(x,n)} = \thDilog(\bar x,-n)^{-1} \qquad &\text{ if } \bpar\in\Real_{>0}\text{.}
	\end{align}
\end{lemma}

\begin{remark}\label{rm:dilogUnit}
	One can see that 
	\begin{align}
	\thDilog(x, -n) = \mathrm{D}_{\bpar^{-1}} (x,n)
	\end{align}
	just by  the Definition \ref{df:thDilog} for $\mathrm{D}_{\bpar^{-1}}$ and carefully substituting $j+n\mapsto  j^\p$. In particular the unitarity for $\bpar>0$ can be re-expressed as
	\begin{align}
	\conj{\thDilog(x,n)} = \left(\mathrm{D}_{\bpar^{-1}} (x,n)\right)^{-1}
	\end{align}  
\end{remark}

\begin{lemma}[Faddeev's difference equations]
	
	Let 
	\begin{equation}\label{eq:chidef}
	\chi^{\pm}(x,n) \equiv e^{2\pi  \frac{\bpar^{\pm 1}}{\sqrt{N}}x} e^{\pm\frac{2\pi i n}{N}}\text{,}
	\end{equation}
	for every $x$, $\bpar\in\Complex$, $\im(\bpar) \neq 0$ $n$, $N\in\Integer$ we have
	\begin{align}
	\thDilog\left(x+i\frac{\bpar^{\pm 1}}{\sqrt{N}}, n\pm 1\right) &= \thDilog\left(x,n\right)\left(1+\chi^{\pm}(x,n)e^{-\pi i \frac{N-1}{N}} e^{\pi i\frac{ \bpar^{\pm2}}{N}}\right)^{-1}\\
	\thDilog\left(x-i\frac{\bpar^{\pm1}}{\sqrt{N}}, n\mp1\right) &= \thDilog(x,n)\left(1+\chi^{\pm}(x,n)e^{\pi i \frac{N-1}{N}} e^{-\pi i\frac{ \bpar^{\pm2}}{N}}\right)
	\end{align}
\end{lemma}
\begin{proposition}\label{pr:ANPoch}
	If $\im(\bpar) >0$ and $\re(\bpar) >0$ we have
	\begin{equation}
	\thDilog(x,n) = \frac{\qPoch{\chi^+(x+\frac{\cth}{\sqrt N},n) }{q^2\omega}}{\qPoch{\chi^-(x-\frac{\cth}{\sqrt N},n)}{\tilde q^2\conj\omega}}
	\end{equation}
	where $q = e^{ i \pi\frac{\bpar^2}{ N}}$, $\tilde q = e^{-\pi i \frac{\bpar^{-2}}{ N}}$, $\omega = e^{\frac{2\pi i }{N}}$ and $\chi^{\pm}(x,n) = e^{2\pi  \frac{\bpar^{\pm 1}}{\sqrt{N}}x} e^{\pm\frac{2\pi i n}{N}}$.
\end{proposition}

\begin{proposition}\label{DthPoles}
	The quantum dilogarithm $\thDilog(x,n)$, for $\im(\bpar)>0$ has poles \\
	$$\begin{cases}
	x = \frac{\cth}{\sqrt N} + i \frac{\bpar^{-1}}{\sqrt N} l +i\frac{\bpar}{\sqrt N} m\\
	n = m - l \mod N 
	\end{cases}$$
	and zeros
	$$\begin{cases}
	x = -\frac{\cth}{\sqrt N} - i \frac{\bpar^{-1}}{\sqrt N} l -i\frac{\bpar}{\sqrt N} m\\
	n = l-m \mod N 
	\end{cases}$$
	for $l$,$m\in\Integer_{>0}$.
	Moreover its residue at $(x_{l,m}, n_{l,m}) = \left(\frac{\cth}{\sqrt N} + i \frac{\bpar^{-1}}{\sqrt N} l +i\frac{\bpar}{\sqrt N} m, m-l\right)$ is 
	\begin{equation}
	\frac{\sqrt{N}}{2\pi\bpar^{-1}} 
	\frac{\qPoch{q^2\omega}{q^2\omega}}{\qPoch{\tilde q^2\conj\omega}{\tilde q^2\conj\omega}}
	\frac{(- \tilde q^2\conj\omega)^l(\tilde q^2 \conj\omega)^{l(l-1)/2}}
	{\Poch{q^2\omega}{q^2\omega}{m}\Poch{\tilde q^2\conj\omega}{\tilde q^2\conj\omega}{l}}
	\end{equation}
\end{proposition}

The following Summation Formula can be shown by a residue computation. It  is well known for $N=1$, i.e. for $\bdilog$, see \cite{FKV1} for example. 
\begin{theorem}[Summation Formula]\label{SummationThm}
	Suppose $\im(\bpar) >0$ and $N$ odd, and let $u$, $v$, $w\in\Complex$ and $a,b,c\in\IntN$ satisfy
	\begin{equation}\label{SummationCond}
	\im\left(v+\frac{\cth}{\sqrt N}\right) >0\text{, } \quad \im\left(-u+\frac{\cth}{\sqrt N} \right) > 0\text{, } \quad\im(v-u) <\im(w) <0\text.
	\end{equation}
	Define 
	\begin{equation}
	\Psi(u,v,w,a,b,c) \equiv \int_{\Affine_N} \frac{\thDilog(x+u, a+d)}{\thDilog (x+v,b+d)}e^{2\pi i wx} e^{-2\pi i \frac{cd}{N}}\diff (x,d)
	\end{equation}
	Then we have that
	\begin{align*}
	\Psi(u,v,w,a,b,c)&\\
	=&\zeta_0 \frac{\thDilog\left(v-u-w+\frac{\cth}{\sqrt N}, b-a-c\right)}{\thDilog\left(-w-\frac{\cth}{\sqrt N}, -c\right)\thDilog\left(v-u+\frac{\cth}{\sqrt N}, b-a\right)} e^{2\pi iw\left(\frac{\cth}{\sqrt N} - u\right)}
	\omega^{ac}\\
	=&\zeta_0^{-1} \frac{\thDilog\left(w+\frac{\cth}{\sqrt N}, c\right)\thDilog\left(-v+u-\frac{\cth}{\sqrt N}, -b+a\right)} {\thDilog\left(-v+u+w-\frac{\cth}{\sqrt N}, -b+a+c\right)}
	e^{2\pi iw\left(-\frac{\cth}{\sqrt N} - v\right)}\omega^{bc}
	\end{align*}
	where $
	\zeta_0 = e^{-\pi i (N - 4\cth^2N^{-1})/12}$.
\end{theorem}
\begin{remark}
	Assumptions (\ref{SummationCond}) even though sufficient are not optimal. Indeed they guarantee the theorem to hold true when the integration is performed along the real line, however we can deform the integration contour as long as 
	\begin{equation}\label{eq:SummationCondRel}
	\abs{\arg (i z)} < \pi -\arg{{\bpar}} \qquad z\text{ being one of }\left\{w,v-u-w, u-v-2\frac{\cth}{\sqrt{N}} \right\}
	\end{equation}
\end{remark}
Using the notation for the Fourier Kernels from \eqref{eq:FourKern} in Appendix \ref{ap:AffineN} we have that
\begin{proposition}[Fourier Transformation Formula, \cite{AK3}]\label{pr:thDilogFourier}
	For $N$ odd we have that
	\begin{align*}
	\int_{\Affine_N} &\thDilog(x,n)\bkt{(x,n);(w,c)}\diff (x,n) =
	\frac{e^{2\pi i w\frac{\cth}{\sqrt N}}}{\thDilog\left(-w-\frac{\cth}{\sqrt N}, -k\right)} e^{-\pi i(N-4\cth^2N^{-1})/12}\\
	&=\thDilog\left(w+\frac{\cth}{\sqrt N}, c\right)\conj{\bkt{(w,c)}} e^{\pi i(N-4\cth^2N^{-1})/12}\\
	\int_{\Affine_N}& (\thDilog(x,n))^{-1}\bkt{(x,n);(w,c)}\diff (x,n) =
	\frac{\bkt{(w,c)}}{\thDilog\left(-w-\frac{\cth}{\sqrt N}, -k\right)} e^{-\pi i(N-4\cth^2N^{-1})/12}\\
	&=\thDilog\left(w+\frac{\cth}{\sqrt N}, c\right)e^{-2\pi i w\frac{\cth}{\sqrt N}} e^{\pi i(N-4\cth^2N^{-1})/12}
	\end{align*}
\end{proposition}
\begin{proposition}[Integral Pentagon Relation]\label{pr:IntPenta}
	Let $\widetilde{\thDilog}(x,n) \equiv \opF_N\circ\mathcal F^{-1}(\thDilog)(x,n)$. We have the following integral relation
	\begin{align*}
	&\bkt{(x,n); (y,m)}\widetilde{\thDilog}(x,n) \widetilde{\thDilog}(y,m)\\
	&\quad=\int_{\Affine_N}\widetilde{\thDilog}(x-z,n-k)\widetilde{\thDilog}(z,k) \widetilde{\thDilog}(y-z, m-k)\bkt{(z,k)}\diff(z,k).
	\end{align*}
\end{proposition}
Before we look at the asymptotic behavior of $\bdilog$ let us recall the classical dilogarithm function, defined on $\abs{z}< 1$ by 
\begin{align}
\Dilog(z) = \sum_{n\geq 1} \frac{z^n}{n^2}
\end{align}
and recall that it admits analytic continuation to $\Complex\setminus[1,\infty]$ through the following integral formula
\begin{align}
\Dilog(z)= - \int_0^z\frac{\log(1-u)}{u}\diff u.
\end{align}

\begin{proposition}
	We have the following behaviour when $\bpar > 0$, $\bpar\rightarrow 0$  and $x$, $n$, $N$ are fixed
	\begin{equation}\label{eq:AsymFaddeev}
	\thDilog(\frac{x}{2\pi \bpar}, n) = \operatorname{Exp}\left(\frac{\Dilog(-e^{x\sqrt{N}})}{2\pi i b^2 N} \right) \phi_x(n)(1+\mathcal{O}(b^2))
	\end{equation}
	where $\phi_x(n)$ is defined by
	$\begin{cases}
	\phi_x(n+1) = \phi_x(n) \frac{(1- e^{x/\sqrt N} \conj\omega^{n+\frac12})}{(1 + e^{x\sqrt N})^{1/N}}\\
	\phi_x(0) = (1+e^{x\sqrt N})^{-\frac{N-1}{2N}} \prod_{j=0}^{N-1} (1 - e^{xN^{-\frac12}} \conj\omega^{j+\frac12})^{\frac jN}
	\end{cases}$
	whenever $N$ is odd.
\end{proposition}
\begin{remark}
	The function $\phi_x$ on the finite set $\IntN$ is a cyclic quantum dilogarithm \cite{FKqdilog},\cite{KQuantTeich}, \cite{K6j}. Precisely $\frac{1}{\phi_x}$ corresponds to the  function $\Psi_{\lambda}$ from Proposition $10$ in \cite{KQuantTeich} with $\lambda = e^{x/\sqrt N}$.
\end{remark}

The Hilbert space $\SqInt (\Affine_N)$ is naturally isomorphic to the tensor product $\SqInt(\Real)\otimes \SqInt(\IntN) \cong \SqInt(\Real) \otimes \Complex^N$, see Appendix \ref{ap:AffineN}.  Let $\opp$ and $\opq$ two self-adjoint operators on $\SqInt(\Real)$ satisfying 
\begin{equation}\label{eq:pqrel}
[\opp,\opq] = \frac{1}{2\pi i}
\end{equation}
and let $X$ and $Y$ unitary operators satisfying 
\begin{equation}\label{eq:finiteoperators}
YX = e^{2\pi i /N} XY \text{, \hskip 1em} X^N = Y^N = 1\text{,}
\end{equation}
together with the cross relations 
\begin{equation}
[\opp,X]=[\opp,Y]=[\opq,X]=[\opq,Y] =0\text{.}
\end{equation}
The equations in \eqref{eq:finiteoperators} imply that $X$ and $Y$ will have finite and the same spectrum, and this will be a subset of the set $\Torus_N$ of all $N$-th complex roots of unity. Let 
$$\LN : \Torus_N \longrightarrow \IntN $$
be the natural group isomorphism. We can define $\LN (A)$, by the spectral theorem, for any operator $A$ of order $N$, such that it formally satisfies 
$$A = e^{2\pi i \LN (A)/N}\text{.}$$
One has that
\begin{equation}\label{eq:sumLNop}
\LN(-e^{-\pi i/N}XY) = \LN(X)+\LN(Y).
\end{equation}
For any function $f:\Affine_N \longrightarrow \Complex$ recall the definition of $\tilde f$ and the operator function $\slashed{f}(\opx,A) \equiv f(\opx,\LN(A))$ from Appendix \ref{ap:AffineN}.
The following Pentagon Identity for $\thDilog$ was first proved in \cite{AK3}, where a projective ambiguity was undetermined and $\abs{\bpar}=1$.
\begin{lemma}[Pentagon Equation]\label{lm:thDilogPent}
	Let $\opp$,$\opq$, $X$ and $Y$ be as above, then the following five-term relation holds
	\begin{equation}\label{eq:thpentagon}
	\SthDilog(\opp,X)\SthDilog(\opq,Y) = \SthDilog(\opq,Y)\SthDilog(\opp+\opq,-e^{\pi i/N}XY)\SthDilog(\opp,X)\text{.}
	\end{equation}
\end{lemma}
\proof
This is equivalent to the Integral Pentagon equation of Proposition \ref{pr:IntPenta}. To see this we need to use equation (\ref{eq:FourierSpectra}) an all the five terms. Then we compare the coefficients of $e^{2\pi iy\opq}Y^{m} e^{2\pi i x\opp}X^{-n}$ and get exactly the integral pentagon equation.\\
An alternative proof follows from the q-Pochammer presentation of $\thDilog$ from Proposition \ref{pr:ANPoch}. 
\endproof

\subsection{Charges}
We are going to define a \emph{charged} version of the dilogarithm. The charges will assume geometrical meaning in the construction of the partition function, however they already satisfy the purpose of turning all the conditional convergent integral relations of the dilogarithm $\thDilog$ (e.g. Proposition \ref{pr:IntPenta} and \ref{pr:thDilogFourier}) into absolutely convergent integrals.\\
Let $a$, $b$ and $c$ be three real positive numbers such that $a+b+c = \frac{1}{\sqrt N}$. We define the charged quantum dilogarithm
\begin{equation}\label{eq:defCharDilog}
\psi_{ a, c}(x,n) := \frac{e^{-2 \pi i \cth ax}}{ \thDilog(x-\cth (a+ c),n)}.
\end{equation}
From the Fourier transformation formula, Proposition \ref{pr:thDilogFourier}, and the inversion formula in Lemma \ref{lm:thDilogInv}, we can deduce the following transformation formulas for $\psi_{ a,c}$ (recall notation \eqref{eq:tildeFourier} for the inverse Fourier transform)
\begin{lemma}\label{lm:dilogSym}
	Suppose $\im(\bpar)(1-\abs{\bpar}) = 0$, then
	\begin{align}
	\tilde{\psi}_{ a,  c}(x,k) &=\psi_{c,  b}(x,k)\bkt{x,k} e^{-\pi i \cth^2 a(a+2c) }\zeta_0 \\
	\conj{\psi_{a, c}(x,k)} &= \psi_{c,  a}(-x,\epsilon k) \bkt{x,k} e^{\pi \cth^2(a+c)^2} \zeta_{N,\text{inv}}\\
	\conj{\tilde{\psi}_{ a,  c}(x, k)} &= \psi_{ b,  c}(-x,\epsilon k)e^{-2\pi i \cth^2ab} \zeta_0
	\end{align}
	where $\zeta_0 = e^{-\pi i(N-4\cth^2N^{-1})/12}$ and $\zeta_{N,\text{inv}} = \zeta_0^{2} e^{-\pi i \cth^2/N}$ and $\epsilon = +1$ if $\bpar>0$ or $\epsilon = -1$ if $\abs{\bpar} = 1$.
\end{lemma}
\begin{remark}
	The hypothesis on positivity of $a$, $b$ and $c$ assure that the Fourier integral of $\tilde{\psi}_{a,c}$ is absolutely convergent.
\end{remark}

\begin{theorem}[Charged Pentagon Equation]
	Let $a_j$, $c_j >0$ such that $\frac{1}{\sqrt{N}} - a_j - c_j >0$ for $j=0$, $1$, $2$, $3$ or $4$. Define $\psi_j\equiv\psi_{ a_j\text{,} b_j} $. Suppose the following relations hold true
	\begin{align}\label{eq:pentaChargesCond}
	a_1 = a_0 +a_2\qquad a_3 = a_2 + a_4 \qquad c_1 = c_0 +a_4 \qquad c_3 = a_0 + c_4 \qquad c_2 = c_1 + c_3\text{}
	\end{align}
	and consider the operators on $\SqInt(\Affine_N)$ defined to satisfy (\ref{eq:pqrel} -\ref{eq:finiteoperators}).
	We have the following charged pentagon relation
	\begin{align}\label{eq:charPenta}
	\psi_1(\opq,\LN(X)) &\psi_3(\opp,\LN(Y)) \xi(a,c) =\\
	\nonumber
	& = \psi_4(\opp,\LN(X)) \psi_2(\opp+\opq, \LN(X)+\LN(Y))\psi_0(\opq,\LN(Y))
	\end{align}
	where $\xi(a,c) =e^{2\pi i \cth^2 (a_0a_2 +a_0 a_4 +a_2 a_4)} e^{\pi i \cth^2a_2^2}$.
\end{theorem}

\section{Quantum Teichm\"uller Theory}\label{sc:QuantumTeich}

In this section we are going to quantize the space $\RatioWell$ following \cite{KQuantTeich}. 
For any fixed $\tau$ the quantization of $\Ratio$ is just the canonical quantization in exponential coordinates of the space $\Real^{M}_{>0}\times\Real^{M}_{>0}$, where $M =  2g -2 + s$, with symplectic form $\omega_{\tau} = \sum_{j=1}^{M} \diff \log u_j\wedge\diff \log v_j$. Formally, following the expectations from canonical quantization of $\Real^{2M}$, we can quantize $\Ratio$ and associate to it an algebra of operator 
\begin{align*}
\mathcal X(\tau) &\text{ generated by }\left\{\hat{u}_j,\, \hat{v}_j  \right\} \text{, where }0\leq j < M\text{,  subject to the relations }\\
&\qquad\hat{u}_j\hat{v}_l = q^{\delta(j-l)} \hat{v}_l\hat{u}_j\quad \hat{u}_j\hat{u}_l = \hat{u}_l\hat{u}_j\quad \hat{v}_j\hat{v}_l = \hat{v}_l\hat{v}_j
\end{align*}
where $q\in\Complex^*$. The algebra we mean here is the associative algebra of non commutative fractions of non commutative polynomials generated by these generators.\\
In order to obtain a quantization of $\RatioWell$ (i.e. triangulation independent) we have to look at the action of the $\PtGr$ generators on coordinates and translate it into an action on the algebras $\mathcal X(\tau)$. Precisely consider the set of the couples $(\tau, \mathcal X(\tau))$ and let the generators  $[\tau, \tau^{\sigma}]$, $[\tau, \rho_i\tau]$ and $[\tau, \omega_{i,j}\tau]$ act on them. The action on the algebras is as follows. The elements $[\tau, \tau^{\sigma}]$ just permutes the indexes of the generators according to the permutation $\sigma$. The change of decoration  $[\tau, \rho_i\tau]$ acts trivially on the operators $(\hat u_j,\hat v_j)$ such that $j\neq i$ and as follows on $(\hat u_i,\hat v_i)$
\begin{align}\label{eq:RotDef}
(\hat u_i,\hat v_i) \mapsto  (q^{-1/2} \hat v_i\hat u_i^{-1},  \hat u_i^{-1})\text{.}
\end{align}
The most interesting generator $[\tau, \omega_{i,j}\tau]$, is again trivial in the triangles not involved in the diagonal exchange, but it maps the two couples of operators $(\hat u_i,\hat v_i)$ and $(\hat u_j,\hat v_j)$ to the two new couples (following formulas \eqref{eq:RatioPtolemy})
\begin{align}\label{eq:tetraEqDef1}
&(\hat u_i,\hat v_i)\bullet (\hat u_j, \hat v_j) \equiv (\hat u_i\hat u_j, \hat u_i\hat v_j + \hat v_i)\\
\label{eq:tetraEqDef2}
&(\hat u_i,\hat v_i)* (\hat u_j, \hat v_j) \equiv ( \hat u_j\hat v_i(\hat u_i\hat v_j+ \hat v_i)^{-1}, \hat v_j(\hat u_i\hat v_j+ \hat v_i)^{-1})\text{.}
\end{align}
 In order to get an actual quantization we need to provide a representation of $\mathcal X(\tau)$ by operators acting on some vector space $\Hilb$. In the original paper \cite{KQuantTeich}, Kashaev proposed representations on the vector spaces $\SqInt(\Real)$ and $\SqInt(\IntN)\simeq \Complex^N$ for $N$ odd. The former was used to construct the Andersen-Kashaev invariants in \cite{AK1}, while the latter are related to the colored Jones polynomials (\cite{K6j}, \cite{MM}) and the Volume Conjecture \cite{KVolKnot}.  In the more recent work \cite{AK3} a representation on the vector space $\SqInt(\Affine_N) \equiv\SqInt(\Real\times\IntN)\simeq \SqInt(\Real)\otimes\Complex^N$ was implicitly proposed, or at least all the basics elements to construct it were presented.
Here we  describe the representations in $\SqInt(\Affine_N)$. 

\subsection{$\SqInt(\Affine_N)$ Representations}\label{sc:ANreps}
Fix $N$ positive odd integer, $\omega \equiv e^{\frac{2\pi i}{N}}$  and $\bpar \in \Complex^*$, $\re(\bpar)>0$.
To each decorated ideal triangle $\overline{\tau}_j\in\tau$ we associate the Hilbert space $ \SqInt(\Affine_N)$. Then the Hilbert space associated to $\Ratio$ will be $\Hilb = \SqInt(\Affine_N)^{\otimes M} \cong \SqInt(\Affine_N^M)$ where $M = 2g -2 + s$ is the number of triangles in $\tau$.
For conventions and notation on the space $\SqInt(\Affine_N)$ see Appendix \ref{ap:AffineN}. For every $i= 0,\dots, M$ let $\opp_i$, $\opq_i$ be self adjoint operator in $\SqInt(\Real)$ and $X_i$, $Y_i$ unitary operators in $\SqInt(\IntN) \simeq \Complex^{N}$ such that
\begin{equation}\label{eq:affineOperator}
[p_i,q_j] = \frac{\delta_{ij}}{2\pi i}\text{, } \qquad Y_iX_j = \omega^{\delta_{ij}} X_jY_i\text{, }\qquad X_i^N = Y_i^N = 1\text{.}
\end{equation}
We can define the operators 
\begin{align}\label{eq:TeichOP2}
&\opu_i = e^{2\pi \frac{\bpar}{\sqrt{N}}\opq_i}Y_i
&\opu_i^* = e^{2\pi \frac{\bpar^{-1}}{\sqrt{N}}\opq_i}Y_i^{-1}\\
&\opv_i = e^{2\pi \frac{\bpar}{\sqrt N}\opp_i}X_i
&\opv_i^* = e^{2\pi \frac{\bpar^{-1}}{\sqrt N}\opp_i}X_i^{-1}
\end{align}
satisfying 
\begin{align}
&\opu_i\opv_j =  q^{\delta_{ij}} \opv_j\opu_i
&\opu_i^*\opv_j^* = \tilde q^{\delta_{ij}}\opv_j^*\opu_i^*\\
&\opu_i\opv_j^* = \opv_j^*\opu_i
&\opu_i^*\opv_j = \opv_j\opu_i^*\\
\label{eq:qdef}
& q = e^{2\pi i \frac{\bpar^2}{N}}\omega
&\tilde q = e^{2\pi i \frac{\bpar^{-2}}{N}}\omega.
\end{align}
The Quantum algebra $\mathcal X(\tau)$  is generated by the $\opu_j$, $\opv_j$ for $j=0,\dots M$, and has a $*$--algebra structure when extended to include $\opu_j^*$ and $\opv_j^*$.
We remark that the $*$ operator we are using here is the standard hermitian conjugation only if $\abs{\bpar} = 1$.\\
Explicitly let $X_j$, $Y_j$, $\opp_j$, $\opq_j$ , $j =1,2$  be operators acting on $\Hilb := \SqInt(\Affine_N^2)$ as follow
\begin{align}\label{eq:TeichOP3}
&\opp_j f(\opx,\opm) = \frac1{2\pi i} \dfield{x_j} f(\opx, \opm)\text{,}\qquad\quad 
\opq_j f(\opx,\opm) = x_jf(\opx,\opm)\\ \label{eq:TeichOP4}
&X_1 f(\opx,\opm) = f(\opx, (m_1 + 1, m_2))\text{,}\quad
X_2 f(\opx,\opm) = f(\opx, (m_1, m_2 +1))\\ \label{eq:TeichOP5}
&Y_j f(\opx, \opm) = \conj{\omega}^{m_j} f(\opx ,\opm)\text{,}
\end{align}
where $\opm = (m_1,m_2)\in \Integer_N^2$ and $\opx = (x_1, x_2) \in \Real^2$.\\
These operators satisfy conditions \eqref{eq:affineOperator}. 
Let $\thPsylog(x,n)\equiv \frac{1}{\thDilog(x,n)}$  and consider the operators
\begin{align}
\opD_{12} &\equiv e^{2\pi i \opq_2\opp_1}\sum_{j,k=0}^{N-1}\omega^{jk}Y_2^jX_1^k\\
\opPsylog_{12} &\equiv \SthPsylog(\opq_1+\opp_2-\opq_2, -e^{-\frac{\pi i}{N}}Y_1X_2\conj{Y_2})\\
\opT_{12}  &\equiv \opD_{12}\opPsylog_{12}     \label{eq:TetraAdef}
\end{align}
One has
\begin{lemma}[Tetrahedral Equations]\label{lm:Teq}
	\begin{align}\label{tetraEq}
	&\opT_{12}\opu_1 = \opu_1\opu_2\opT_{12} \\
	&\opT_{12}\opv_1\opv_2 = \opv_2\opT_{12} \\
	&\opT_{12}\opv_1\opu_2 = \opv_1\opu_2\opT_{12}\\ 
	&\opT_{12}\opv_1 = (\opu_1\opv_2+\opv_1) \opT_{12}\\
	&\opT_{12}\opT_{13}\opT_{23} = \opT_{23}\opT_{12} \label{eq:TetraPenta}
	\end{align}
\end{lemma}
\begin{remark}
	If we define $\tilde{\opT}_{12} = \opD_{12} \tilde{\opPsylog}_{12}$ where
	\begin{equation}
	\tilde{\opPsylog} \equiv \slashed{\Psi}_{\bpar^{-1}}(\opq_1+\opp_2-\opq_2, -e^{-\frac{\pi i}{N}}\conj{Y}_1\conj{X}_2{Y_2})
	\end{equation}
	then $\tilde\opT$ satisfies equations (\ref{tetraEq} -- \ref{eq:TetraPenta}) with $\opu_i$ and $\opv_i$ substituted by $\opu^*_i$ and $\opv^*_i$.
	However from Remark \ref{rm:dilogUnit} we know that $\Psi_{\bpar^{-1}} (x, n) = \Psi_{\bpar}(x, -n)$, and 
	$$\left(-e^{-\frac{\pi i}{N}}Y_1X_2\conj{Y_2}\right)^{-1} = -e^{\frac{\pi i}{N}}\conj Y_1 Y_2\conj X_2 =-e^{-\frac{\pi i}{N}}\conj{Y}_1\conj{X}_2{Y_2}
	$$
	so that 
	$$\tilde{\opT} = \opT.$$
\end{remark}

From Lemma \ref{lm:Teq} we have the following implementations of equations (\ref{eq:RotDef} -- \ref{eq:tetraEqDef2}).
\begin{proposition}
	Let $\opw_i\equiv (\opu_i,\opv_i)$ and $\opw^*_i =  (\opu_i^*,\opv_i^*)$. Then we have 
	\begin{align}
	&\opw_1\bullet\opw_2 \opT_{12} = \opT_{12} \opw_1\text{,\hskip 1em} \opw_1*\opw_2\opT_{12} = \opT_{12}\opw_2\text{,}\\
	&\opw_1^*\bullet\opw_2^* \opT_{12} = \opT_{12} \opw_1^*\text{,\hskip 1em} \opw_1^**\opw_2^*\opT_{12} = \opT_{12}\opw_2^*\text{.}
	\end{align}
\end{proposition}

\begin{proposition}
	Let
	\begin{align}
	&\opA\equiv  e^{3\pi i \opq^2} e^{\pi i (\opp+ \opq)^2} \sum_{j=0}^{N-1} \conj{\bkt{j}}^3 Y^{3j}\sum_{l= 0}^{N-1} \conj{\bkt{l}} (-e^{-\pi i /N }YX)^l,
	\end{align}
	where $\bkt{n} = e^{-\pi i n(n+N)/N}$ and $Y$ and $X$ are as above.
	Then 
	\begin{align}
	&\opA(\opu, \opv) = (q^{-1/2} \opv\opu^{-1}, \, \opu^{-1})
	&\opA(\opu^*, \opv^*) = (\tilde q^{-1/2} \opv^*(\opu^*)^{-1}, \, (\opu^*)^{-1})
	\end{align}
	where $q$ and $\tilde q$ are defined by equation \eqref{eq:qdef}.
\end{proposition}

\section{Quantization of the Model Space for Complex Chern-Simons Theory}\label{sc:ModelSpace}
In this Section we want to quantize the space  
$\Complex^*\times\Complex^*$ with the complex symplectic form 
$$\omega_{\Complex} = \frac{\diff x\wedge\diff y}{xy}\text{.}$$ 
We think of it as a model space for Complex Chern-Simons Theory because it is an open dense of the $\PSL(2,\Complex)$ moduli space of flat connections on a four punctured sphere, with unipotent holonomy around the punctures, \cite{AK3,KPSL2R,D3D,FGcluster}. 
Tetrahedral operators are supposedly related to states in the quantization of the four punctured sphere. Since we want to construct knot invariants starting from tetrahedral ideal triangulations this is the space we need to quantize. We follow the ideas in Andersen and Kashaev \cite{AK3} using a real polarization with contractible leaves. 
We will further show that the level $N$  Weil-Gel'fand-Zak Transform relates this quantization  with the $\SqInt(\Affine_N)$ representations in Quantum Teichm\"uller theory. To use this transform  to relates  the Andersen--Kashaev invariants to complex Chern--Simons Theory was already proposed in \cite{AK3}. However the relation between the two approaches was not as tight as the one present here.

Let $t = N+is \in\Complex^*$ be the quantization constant, for $N\in\Real$ and $s\in \Real\sqcup\,i\Real$. Denote also $\tilde t = N-is$.
Fix $\bpar\in\Complex$ such that $s= -iN\frac{1-\bpar^2}{1+\bpar^2}$  and $\re\bpar >0$. This substitution, for $s\in i\Real$, is only possible when $-N<is<N$. Notice that 
\begin{align}
&s\in\Real\,\, \iff \,\, \abs{\bpar}= 1\text{ and }\bpar\neq \pm i, &s\in\, i\Real \,\, \iff \,\, \im\bpar = 0
\end{align} 
and notice the following useful expressions
\begin{align}
t = \frac{2N}{1+\bpar^2}, \qquad \qquad \tilde t = \frac{2N}{1+\bpar^{-2}}.
\end{align}
Consider the covering maps 
\begin{align}\label{eq:expCoordDef}
\zeta^{\pm} \colon \Real^2& \longrightarrow \Complex^*\\
\nonumber
(z, n)  &\mapsto\exp\left(2\pi \bpar^{\pm 1}z \pm 2\pi i n  \right)
\end{align}
and consequently
\begin{equation}\label{eq:expProjDef}
\pi^{\pm}:\Real^2\times\Real^2 \longrightarrow \Complex^*\times\Complex^*\text{, } \qquad \pi^{\pm} = (\zeta^{\pm},\,\zeta^{\pm})
\end{equation}
such that 
\begin{align}
&\Complex^*\times\Complex^* \,\ni\, (x,y) = \pi^+((z,n),(w,m))\text{, } \\
\nonumber
&\Complex^*\times\Complex^* \,\ni\, (\tilde x,\tilde y) = \pi^-((z,n),(w,m))\quad \text{for } ((z,n),(w,m))\in\Real^2\times\Real^2\text{.}
\end{align}
We remark that 
\begin{equation}
\conj{\zeta^+(z,n)} = \zeta^-(z,n) \qquad \iff \qquad \abs{\bpar} = 1
\end{equation}
in this case $\pi^{-} = \conj{\pi^+}$ and $\tilde x = \conj x$, $\tilde y= \conj y$. In this sense $x$, $y$ $\tilde x$ and $\tilde y$ are natural coordinate functions to quantize in $\Complex^*\times\Complex^*$. If $\bpar\in \Real$ they are still coordinates functions for the underlying real manifold, but we lose the complex conjugate interpretation.
We will first consider the quantization of the covering $\Real^2\times\Real^2$. Define the form 
\begin{align}
\omega_t \equiv \frac{t}{4\pi}(\pi^{+})^*(\omega_{\Complex}) +\frac{\tilde t}{4\pi} (\pi^-)^*(\omega_{\Complex}).
\end{align}  

\begin{lemma}\label{lm:omegatExpl}
	\begin{equation}
	\omega_t = 2\pi N(\diff z \wedge\diff w -\diff n \wedge \diff m)\text{.}
	\end{equation}
	In particular it is a real symplectic $2$ form on $\Real^2\times\Real^2$, independent of $\bpar$.
\end{lemma} 

Over $\Real^2\times\Real^2$ we take the trivial line bundle $\widetilde \Line = \Real^2\times\Real^2\times \Complex$. On the $N$-th tensor power of this line bundle $\widetilde\Line^N$  we consider the connection
\begin{equation}
\nabla^{(t)} \equiv \diff - i\alpha_t
\end{equation}
where 
\begin{align}
&\alpha_t \equiv \frac{t}{4\pi} \alpha_{\Complex}^{+} + \frac{\tilde t}{4\pi} \alpha_{\Complex}^{-} \text{, }\\
&\alpha_{\Complex}^{\pm} \equiv 2\pi^2(\bpar^{\pm 1} z \pm in)\diff (\bpar^{\pm 1}w \pm im) -2\pi^2 (\bpar^{\pm 1}w \pm im)\diff (\bpar^{\pm 1} z \pm in)
\end{align}
In analogy to Lemma \ref{lm:omegatExpl} we have
\begin{align}
\alpha_t = \pi N(z\diff w -w\diff z -n\diff m + m\diff n).
\end{align}
It is simple to see that
\begin{align}
&\diff\alpha_{\Complex}^{\pm}  = (\pi^{\pm})^*(\omega_{\Complex})\text{, }\qquad \text{  which implies}\\
& F_{\nabla^{(t)}} = -i\omega_t\text{.}
\end{align} 
Further, on $\Real^2\times\Real^2$ we have an action of $\Integer\times\Integer$ compatible with the projection $\pi^+$, i.e.
\begin{align}
\left(\Integer\times\Integer\right)\times \left(\Real^2\times\Real^2\right) &\longrightarrow \Real^2\times\Real^2\\
\nonumber
(\lambda_1,\lambda_2)\cdot((z,n), (w,m))&\mapsto ((z,n+ \lambda_1), (w,m+ \lambda_2))\\
\end{align} 
that  satisfies 
\begin{align}
\pi^{\pm}((z,n+ \lambda_1), (w,m+ \lambda_2))& = \pi^{\pm}((z,n), (w,m))
\end{align}
This action can be lifted to an action $\widetilde{\Line}^N$ in such a way that the quotient bundle $\Line^N \equiv\widetilde{\Line}^N\ /(\Integer)^2 \rightarrow \Real^4/\Integer^2$ has first chern class $c_1(\Line^N) = \frac{1}{2\pi} \left[\omega_t\right]$ ($\omega_t$ is evidently $\Integer^2$--invariant). Such a condition gives the requirement  (which is in fact the pre-quantum condition) $\frac{1}{2\pi} \left[\omega_t\right]\in H^2((\Real^2\times\Real^2)/(\Integer^2) \text{, }\Integer)$, which boils down to the requirement $N\in \Integer$. Explicitly the action of $\Integer \times\Integer$ on $\widetilde{\Line}^N$ is given by the following two multipliers
\begin{align}\label{eq:multipliers}
&e_{(1,0)} = e^{-\pi N im}, &e_{(0,1)} = e^{\pi  Nin}.
\end{align} 
That means that we consider the space of sections 
\begin{equation}
(\Smooth(\Real^4, \widetilde \Line^N))^{\Integer^2}
\end{equation}
of $\Integer^2$--invariant, smooth sections of $\widetilde{\Line}^N$ .
Explicitly 
\begin{align}
\nonumber
&s\in(\Smooth(\Real^4, \widetilde \Line^N))^{\Integer^2} \text{  if and only if  } s\in\Smooth(\Real^4, \widetilde \Line^N)
\text{   and satisfies } \\
&s((z,n+1),(w,m)) = e^{-\pi i N m}s((z,n),(w,m))\text{,}\\
& s((z,n),(w,m+1)) = e^{\pi i N n} s((z,n),(w,m))
\end{align}

\begin{lemma}
	\begin{align*}
	\nabla^{(t)}s \,\in\, (\Smooth(\Real^4, \widetilde \Line^N))^{\Integer^2}
	\text{, }\quad \text{ for any } s\,\in\, (\Smooth(\Real^4, \widetilde \Line^N))^{\Integer^2}
	\end{align*}
\end{lemma}

The following Hermitian structure on $\widetilde\Line^N$ is $\Integer^2$--invariant and parallel with respect to $\nabla^{(t)}$.
\begin{equation}
s\cdot s^{\prime} (p) \equiv s(p)\conj{s^{\p}(p)}\text{,  } \quad \text{for any }p\in \Real^2\times\Real^2
\end{equation}
Being parallel here means that 
\begin{equation}
\diff (s\cdot s^{\p}) = (\nabla^{(t)}s)\cdot s^\p + s\cdot (\nabla^{(t)}s^\p)\text{,}
\end{equation}
and this is a simple consequence of $\alpha_t$ being a real $1$-form.
It follows that the following is a well defined inner product in the completion of $\left(\left(\SqInt\cap\Smooth\right)(\Real^4, \widetilde \Line^N)\right)^{\Integer^2}$

\begin{align}\label{eq:innerPreQuant}
\left(s\text{,} s^\p\right) \equiv \int_{\Real}\diff z\int_{\Real }\diff w\left( \int_0^1\diff n\int_0^1\diff m \,\, s\cdot s^{\prime} \right)
\end{align}

\begin{lemma}
	We have the following Hamiltonian vector field for the coordinates functions on $\Real^2\times\Real^2$
	\begin{align*}
	&X_z = \frac{1}{2\pi N }\dfield{w} 
	&X_w = -\frac{1}{2\pi N }\dfield{z}\\
	&X_n = -\frac{1}{2\pi N }\dfield{m}
	&X_m = \frac{1}{2\pi N }\dfield{n}
	\end{align*}
\end{lemma}

From the definition of the pre-quantum operator $\hat f$ associated to the observable $f$, we have
\begin{equation}
\hat{f} = -i\nabla_{X_f} + f
\end{equation}
\begin{lemma}[Pre--Quantum operators]\label{lm:pqOps}
	The following are the pre--quantum operators for the coordinate functions on $\Real^2\times\Real^2$
	\begin{align*}
	&\hat z = \frac{-i}{2\pi N }\nabla^{(t)}_w + z  
	&\hat w = \frac{i}{2\pi N }\nabla^{(t)}_z + w\\
	&\hat n = \frac{i}{2\pi N }\nabla^{(t)}_m + n
	&\hat m = \frac{-i}{2\pi N }\nabla^{(t)}_n + m
	\end{align*}
	and they satisfy the following canonical commutation relations 
	\begin{align}
	&\comm{\hat{z}\text{, } \hat{w}} = \frac{1}{2\pi iN}
	&\comm{\hat n\text{, } \hat m} = -\frac{1}{2\pi iN}\\
	&\comm{\hat z, \hat n} =\comm{\hat z, \hat m} = \comm{\hat w,\hat n} = \comm{\hat w, \hat m} = 0 
	\end{align}
\end{lemma}
The Hermitian line bundle $ \Line^N\rightarrow \left(\Real^2\times\Real^2/\Integer^2\right)$ together with the connection $\nabla^{(t)}$ define a pre--Quantization of the theory. In order to finish the quantization program we need to choose a Lagrangian polarization.

Choose the following real  Lagrangian polarization 
\begin{equation}
\widetilde{\mathcal P} \equiv \operatorname{Span}_{\Real}\left\{\dfield{w}+ \dfield{n}\text{,  } \dfield{z} - \dfield{m} \right\}.
\end{equation}
The leaves of this polarization are all contractible after quotient by the action of $\Integer^2$ on $\Real^2\times\Real^2$, so we do have polarized global sections. In particular the space $T\subset \Real^2 \times\Real^2$ 
\begin{equation}
T\equiv\left\{z = w = 0 \right\}
\end{equation}
is a transversal for the polarization. For any $\psi \in (\Smooth(\Real^4, \widetilde \Line^N))^{\Integer^2}$ polarized by $ \widetilde{\mathcal P}$, the following two differential equations will determine $\psi \equiv \psi ((z,n),(w,m))$ by its value in $(n,m)$
\begin{align}
&\nabla^{(t)}_w\psi = -\nabla_n^{(t)}\psi   &\nabla_z^{(t)}\psi = \nabla^{(t)}_m\psi\text{.}
\end{align} 
The space $T/\Integer^2$ is of course $\Torus\times\Torus$, and the line bundle $\Line^{N}$ will restrict to a non trivial line bundle over $\Torus\times\Torus$ that we shall call $\Line^N$ again. The quantum space that we obtain is then 
\begin{align}
{\Hilb}^{(N)} \equiv \Smooth\left(\Torus\times\Torus, \Line^N\right).
\end{align} 
We consider the obvious inner product on $\Hilb^{(N)}$
\begin{align}
\left(\psi\text{, }\phi\right) = \int_0^1\int_0^1  \psi\conj{\phi}\,\,\diff n\diff m
\end{align}
that is  the standard inner product on the completion $\SqInt\left(\Torus\times\Torus, \Line^N\right)$.
Finally the quantum operators acts on polarized sections as  
\begin{align}
& \hat x\equiv\exp\left( 2\pi \bpar \hat z + 2\pi i\hat n \right)
=  \exp\left(i\frac{\bpar}{N} \nabla^{(t)}_n - \frac{1}{N}\nabla_m^{(t)} + 2\pi i n\right) \\
& \hat y\equiv\exp\left( 2\pi \bpar\hat w + 2\pi i\hat m \right)
= \exp\left( i\frac{\bpar}{N} \nabla^{(t)}_m + \frac{1}{N}\nabla_n^{(t)} +  2\pi im \right)\\
&\hat {\tilde x}\equiv \exp\left(2\pi\bpar^{-1}\hat z -2\pi i\hat n \right)
=\exp\left(  i\frac{\bpar^{-1}}{N} \nabla^{(t)}_n + \frac{1}{N}\nabla_m^{(t)} - 2\pi in \right)\\
&{\hat {\tilde y }}\equiv\exp\left( 2\pi\bpar^{-1}\hat w -2\pi i\hat m\right)
= \exp\left( i\frac{\bpar^{-1}}{N} \nabla^{(t)}_z - \frac{1}{N}\nabla_n^{(t)}  - 2\pi im \right)
\end{align}
Now we are going to connect the quantization with the Quantum Teichm{\"u}ller theory.
Recall the operators $\opu = \opu(\bpar)$ and $\opv = \opv(\bpar)$  from equations (\ref{eq:affineOperator} -- \ref{eq:TeichOP5}),
and recall that they depend on a parameter $\bpar$. Define the rescaling operator
\begin{align}
\nonumber
\mathcal O_{\sqrt N} \colon \SqInt(\Affine_N) &\longrightarrow \SqInt(\Affine_N)\\
\mathcal O_{\sqrt N}(\opf)(z,n) &= \opf\left(\sqrt N z, n\right)
\end{align} 
and the following rescaled analogues of the Quantum Teichm{\"u}ller Theory operators 
\begin{align}
&\hat u = \mathcal O_{\sqrt N} \circ  \opu(\bpar^{-1} ) \circ \mathcal O_{\sqrt N}^{-1}
&\hat v = \mathcal O_{\sqrt N}\circ \opv(\bpar^{-1}) \circ \mathcal O_{\sqrt N}^{-1}\\
&\hat{u}^* = \mathcal O_{\sqrt N} \circ  \opu^*(\bpar^{-1} ) \circ \mathcal O_{\sqrt N}^{-1}
&\hat{v}^* = \mathcal O_{\sqrt N}\circ \opv^*(\bpar^{-1}) \circ \mathcal O_{\sqrt N}^{-1}
\end{align}
which acts on $\opf\in\SqInt(\Affine_N)$ as
\begin{align}
&\hat{u}^* \mathsf f (z,l) = e^{2\pi \bpar z}e^{2\pi il /N}\mathsf f(z,l) 
&\hat u  \mathsf f (z,l) = e^{2\pi \bpar^{-1} z}e^{-2\pi il /N}\mathsf f(z,l)\\
&\hat{v}^* \mathsf f (z,l) = f(z-i\frac{\bpar}{N}, l -1)
&\hat v \mathsf f (z,l) = f(z-i\frac{\bpar^{-1}}{N}, l +1)
\end{align}

We make use of the level-$N$  Weil-Gel'fand-Zak Transform, \cite{AK3}.
\begin{theorem}
	Recall the line bundle $\Line^N$.
	The following map $Z^{(N)} \colon \Schw(\Affine_N)\longrightarrow \Smooth(\Torus\times\Torus, \Line^N)$ is a an isomorphism
	\begin{align}
	Z^{(N)} (\mathsf f) (n,m)= \frac{1}{\sqrt N} e^{\pi i N mn}\sum_{p\in\Integer}\sum_{l=0}^{N-1}\mathsf f\left(n+\frac{p}{N},l\right)e^{2\pi i mp} e^{2\pi i lp/N}
	\end{align}
	with inverse
	\begin{align}
	\nonumber
	\conj{Z}^{(N)} (s) (x, j) &= \frac{1}{\sqrt N} \sum_{l=0}^{N-1} e^{-2\pi i \frac{lj}{N} } \int_0^1 s\left(x - \frac lN, v\right) e^{-\pi i N(x +\frac lN )v} \diff v\text{.}
	\end{align}
	which preserves the inner products  $\SqInt(\Affine_N)$ and $\left(\cdot,\cdot\right)$, i.e.
	\begin{align*}
	\left(Z^{(N)} (\mathsf f)\text{, }Z^{(N)} (\mathsf g)\right) = \bkt{ \mathsf f\text{, } \mathsf g}
	\end{align*}
	and so extends to an isometry between $\SqInt(\Affine_N)$ and $\SqInt(\Torus\times\Torus, \Line^N)$.
\end{theorem}

\begin{proposition}\label{pr:EquivOpers}
	We have
	\begin{align*}
	&Z^{(N)}\circ \hat{u}^* \circ  (Z^{(N)})^{-1} = \hat y^{-1} 
	&Z^{(N)}\circ \hat{v}^*\hat \circ  (Z^{(N)})^{-1} = \hat x^{-1}\\
	&Z^{(N)}\circ \hat{u}  \circ (Z^{(N)})^{-1} = \hat{\tilde y} ^{-1}  
	&Z^{(N)}\circ \hat{v}  \circ  (Z^{(N)})^{-1} = {\hat{\tilde x} }^{-1}
	\end{align*}
\end{proposition}

All together we have showed that the quantization for the model space of complex Chern-Simons theory is equivalent to the $\SqInt(\Affine_N)$ representations of the quantum algebra defined from Quantum Teichm\"uller Theory. 
In the following section we will extend the $\SqInt(\Affine_N)$ representations to knots invariants following the recipe given by Andersen and Kashaev in \cite{AK1}. The previous discussion on the different quantizations serves to link such invariants to Complex Quantum Chern--Simons Theory.

\section{The Andersen--Kashaev Teichm\"uller TQFT at Level $N$}\label{ch:TQFT}

\subsection{Angle Structures on $3$-Manifolds}
We present here shaped triangulated pseudo $3$-manifolds, which are the combinatorial data underlying the Andersen-Kashaev construction of their invariant. Following strictly \cite{AK1} we will describe the \emph{categoroid} of \emph{admissible} oriented triangulated pseudo $3$-manifolds, where the words admissible and categoroid go together because admissibility is what will obstruct us to have a full category. See Appendix \ref{ap:categoroid} for a definition of categoroids.

\begin{definition}[Oriented Triangulated Pseudo $3$-manifold]
	An \emph{Oriented Triangulated Pseudo $3$-manifold }$X$ is a finite collection of $3$-simplices (tetrahedra) with totally ordered vertices together with a collection of \emph{gluing homeomorphisms} between some pairs of codimension $1$ faces, so that every face is in, at most, one such pairs. By gluing homeomorphism we mean a vertex order preserving, orientation reversing, affine homeomorphism between the two faces.\\
	The quotient space under the glueing homeomorphisms has the structure of CW-complex with oriented edges.
\end{definition}

For $i\in\{0,1,2,3\}$ we denote by $\Delta_i(X)$ the collection of $i$-dimensional simplices in $X$ and, for $i>j$, we denote 
$$\Delta_i^j(X) = \{ (a,b) | a\in \Delta_i(X), b\in\Delta_j(a)\}\text{.}$$
We have projection maps 
$$\phi_{i,j}:\Delta_i^j(X) \longrightarrow \Delta_i(X)\text{, \hskip 1em } \phi^{i,j}:\Delta_i^j(X) \longrightarrow \Delta_j(X)\text{,}$$
and boundary maps
$$\partial_i:\Delta_j(X) \longrightarrow \Delta_{j-1} (X) \text{, \hskip 1em } \partial_i[v_0,\dots , v_j]\mapsto [v_0,\dots,v_{i-1},v_{i+1}, \dots , v_j]$$ 
where $[v_0,\dots, v_j]$ is the $j$-simplex with vertices $v_0,\dots , v_j$ and $i\leq j$.

\begin{definition}[Shape Structure]
	Let $X$ be an oriented triangulated pseudo $3$-manifold. A \emph{Shape Structure} is a map 
	$$\alpha_X : \Delta_3^1(X) \longrightarrow \Real_{>0}\text{,}$$
	so that, in every tetrahedron, the sum of the values of $\alpha_X$ along three incident edges is $\pi$.\\
	The value of the map $\alpha_X$ in an edge $e$ inside a tetrahedron $T$ is called the \emph{dihedral angle} of $T$ at $e$.
	If we allow $\alpha_X$ to take values in $\Real$ we talk about a \emph{Generalized Shape Structure}.\\
	The set of shape structures supported by $X$ is denoted $\Shp (X)$. The space of generalized shape structures is denoted by $\GShp (X)$.
	$X$ together with $\alpha_X$ is called \emph{Shaped Pseudo $3$-manifold}.
\end{definition}

\begin{remark}[Ideal Tetrahedron]\label{rm:idealTetra}
	A shape structure on a simplicial tetrahedron $T$ as above defines an embedding of $T\setminus\Delta_0(T)$ in the hyperbolic $3$--space $\Hyp{3}$ which extends to a map of $T$ to $\conj{\Hyp{3}}$. In fact  we can change a given embedding, so that it send the four vertices $(v_0, v_1, v_2, v_3)$ to the four points $(\infty,0, 1, z )\in\Complex\mathbb P^1 \simeq \partial \Hyp{3}$, where
	\begin{equation*}
	z = \frac{\sin \alpha_T([v_0,v_2])}{\sin \alpha_T([v_0,v_3])} \exp\left(i\alpha_T([v_0,v_1])\right)\text{.}
	\end{equation*}
	This four points in $\partial\Hyp{3}$ extend to a unique ideal tetrahedron in $\Hyp{3}$, by taking the geodesic convex hull, that has dihedral angles defined by $\alpha_T$.
\end{remark}

\begin{remark}
	In every tetrahedron, its orientation induces a cyclic ordering of all triples of edges meeting in a vertex. Such a cyclic ordering descends to a cyclic ordering of the pairs of opposite edges of the whole tetrahedron. Moreover, it follows from the definition that opposite edges share the same dihedral angle. Hence, we get a well defined cyclic order preserving projection $p: \Delta_3^1(X) \longrightarrow \Delta_3^{1/p} (X) $ which identifies opposite edges. $\alpha_X$ descends to a map from $\Delta_3^{1/p} (X) $ and we can consider the following skew-symmetric functions 
	$$\varepsilon_{a,b} \in \{0,\pm1\} \text{, \hskip 1em} \varepsilon_{a,b} = - \varepsilon_{b,a} \text{, \hskip 1em} a,b \in\Delta_3^{1/p} (X)\text{,}$$ 
	defined to be $\varepsilon_{a,b} =0$ if the underlying tetrahedra are distinct, and $\varepsilon_{a,b} = +1$ if the underlying tetrahedra coincides and $b$ cyclically follows $a$ in the order induced on $\Delta_3^{1/p} (X)$.
\end{remark}

\begin{definition}
	To any shaped pseudo $3$-manifold $X$, we associate a \emph{Weight} function
	$$\omega_X: \Delta_1(X) \longrightarrow \Real_{>0} \text{, \hskip 1em } \omega_X(e) = \sum_{a\in(\phi^{3,1})^{-1}(e)} \alpha_X(a)\text{.}$$
	An edge $e$  in $X$ is called \emph{balanced} if $e$ is internal and $\omega_X(e) = 2\pi$. A shape structure is fully balanced if all its edges are balanced.
\end{definition}

The shape structures of closed fully balanced $3$-manifolds are called \emph{Angle Structures} in the literature. For more details on them and their geometric admissibility see \cite{Lack} and \cite{LT}.

\begin{definition}
	A \emph{leveled} (generalised) shaped pseudo $3$-manifold is a pair $(X, l_X)$ consisting of a (generalized) shaped pseudo $3$-manifold $X$ and a real number $l_X\in\Real$, called the \emph{level}.
	The set of all leveled (generalised) shaped pseudo $3$-manifolds is denoted by $\LShp (X)$ (respectively $\LGShp (X)$).
\end{definition}

There is a gauge action of $\Real^{\Delta_1(X)}$ on $\LGShp (X)$.
\begin{definition}\label{df:ggShp1}
	Let $(X, l_X)$ and $(Y, l_Y)$ be two (generalized) leveled shaped pseudo $3$-manifolds. They are said to be \emph{gauge equivalent} if there exists an isomorphism $h:X\longrightarrow Y$  of the underlying cellular structures, and a function $g:\Delta_1(X) \longrightarrow \Real$  such that 
	\begin{align*} 
	&\Delta_1(\partial X) \subset g^{-1}(0)\text{,} \\
	&\alpha_Y (h(a)) = \alpha_X(a) + \pi \sum_{b\in\Delta_3^1(X)} \varepsilon_{p(a),p(b)} g(\phi^{3,1}(b))\text{, \hskip 1em } \forall a\in\Delta_3^1(X)\text{, and}	\\
	&l_Y = l_X + \sum_{e\in\Delta_1(X)} g(e) \sum_{a\in(\phi^{3,1})^{-1} (e)} (\frac{1}{3} - \frac{\alpha_X(a)}{\pi})\text{.}
	\end{align*}
\end{definition}

We remark that $\omega_X = \omega_Y\circ h$.

\begin{definition}
	Let $(\alpha_X, l_X)$ and $(\alpha_{X^\p}, l_{X^\p})$ be two (generalized) leveled  shape structures of the oriented pseudo $3$-manifold $X$. They are said \emph{based gauge equivalent} if they are gauge equivalent as in Definition \ref{df:ggShp1} if the isomorphism $h:X\longrightarrow X$ is the identity.
\end{definition}

Based gauge equivalence is an equivalence relation in the sets $\Shp (X)$, $\LShp (X)$, $\GShp (X)$, $\LGShp (X)$  and the quotient sets are denoted (resp.)  $\Shp_r (X)$, $\LShp_r (X)$, $\GShp_r (X)$, $\LGShp_r (X)$.
We remark that $\Shp_r (X)$ is an open convex (possibly empty) subset of the space $\GShp_r (X)$. We will return to existence of shape structures later. Let us focus on $\GShp(X)$ for now.
Let $$\tilde{\Omega}_X : \GShp (X) \longrightarrow \Real^{\Delta_1(X)}$$ be the map which sends the shape structure $\alpha_X$ to the corresponding weight function $\omega_X$. This map is gauge invariant, so it descends to a map $$\tilde{\Omega}_{X,r} : \GShp_r (X) \longrightarrow \Real^{\Delta_1(X)}$$
For fixed  $a\in\Delta_3^{1/p}(X)$ we can think of $\alpha_a :=\alpha_X(a)$ as an element of $C^{\infty}\left(\GShp(X)\right)$.

\begin{definition}[\cite{NZ}]
	The Neumann-Zagier symplectic structure on $\GShp (X)$ is the unique symplectic structure which induces the Poisson bracket $\{ \cdot,\cdot \}$ satisfying 
	$$\{ \alpha_a,\alpha_b\} = \varepsilon_{a,b}$$
	for all $a$, $b\in\Delta_3^{1/p}(X)$.
\end{definition}

For a triangulated pseudo $3$-manifold we have a symplectic decomposition 
$$\GShp (X) = \prod_{T\in\Delta_3(X)} \GShp(T)\text{.}$$

\begin{theorem}[\cite{AK1}]
	The gauge action of $\Real^{\Delta_1(X)}$ on $\GShp (X)$ is symplectic and $\tilde{\Omega}_X$ is a moment map for this action.
	It follows that $\GShp_r (X)= \GShp (X)/\Real^{\Delta_1(X)}$ is a Poisson manifold with symplectic leaves corresponding to the fibers of  $\tilde{\Omega}_{X,r}$.
\end{theorem}

Let $N_0(X)$ be a sufficiently small neighbourhood of $\Delta_0(X)$, then $\partial N_0(X)$ is a surface which inherits a triangulation from $X$, with a shape structure, if $X$ has a shape structure. Notice that this surface can have boundary if $\partial X \neq \emptyset$.

\begin{theorem}[\cite{AK1}]\label{th:shp1}
	The map $$\tilde{\Omega}_{X,r} : \GShp_r (X)\longrightarrow \Real^{\Delta_1(X)}$$ is an affine $\Hom^1(\partial N_0(X),\Real)$-bundle. The Poisson structure of $\GShp_r(X)$ coincide with the one induced by the $\Hom^1(\partial N_0(X),\Real)$-bundle structure.\\
	If $h:X\longrightarrow Y$ is an isomorphisms of cellular structure, the induced morphism $h^*:\GShp_r(Y)\longrightarrow\GShp_r(X)$ is compatible with all this structures, i.e. it is a Poisson affine bundle morphism which fiberwise coincide with the naturally induced group morphism $h^*:\Hom^1(\partial N_0(Y),\Real)\longrightarrow \Hom^1(\partial N_0(X),\Real)$. Moreover $h^*$ maps $\Shp_r(Y)$ to $\Shp_r(X)$.
\end{theorem}

\begin{definition}[Shaped $3-2$ Pachner moves] 
	Let $X$ be a shaped pseudo $3$ manifold and let $e$ be a balanced internal edge in it, shared exactly by three distinct tetrahedra $t_1$, $t_2$ and $t_3$ with dihedral angles at $e$ exactly $\alpha_1$, $\alpha_2$ and $\alpha_3$.  Then the triangulated pseudo $3$-manifold $X_e$ obtained by removing the edge $e$, and substituting the three tetrahedra $t_1$, $t_2$ and $t_3$ with other two new tetrahedra $t_4$ and $t_5$ glued along one face, is topologically the same space as $X$. In order to have the same weights of $X$  on $X_e$, the dihedral angles of $t_4$ and $t_5$ are uniquely determined by the ones of $t_1$, $t_2$ and $t_3$  as follows
	\begin{figure}[ht!]
		\centering
		\includegraphics{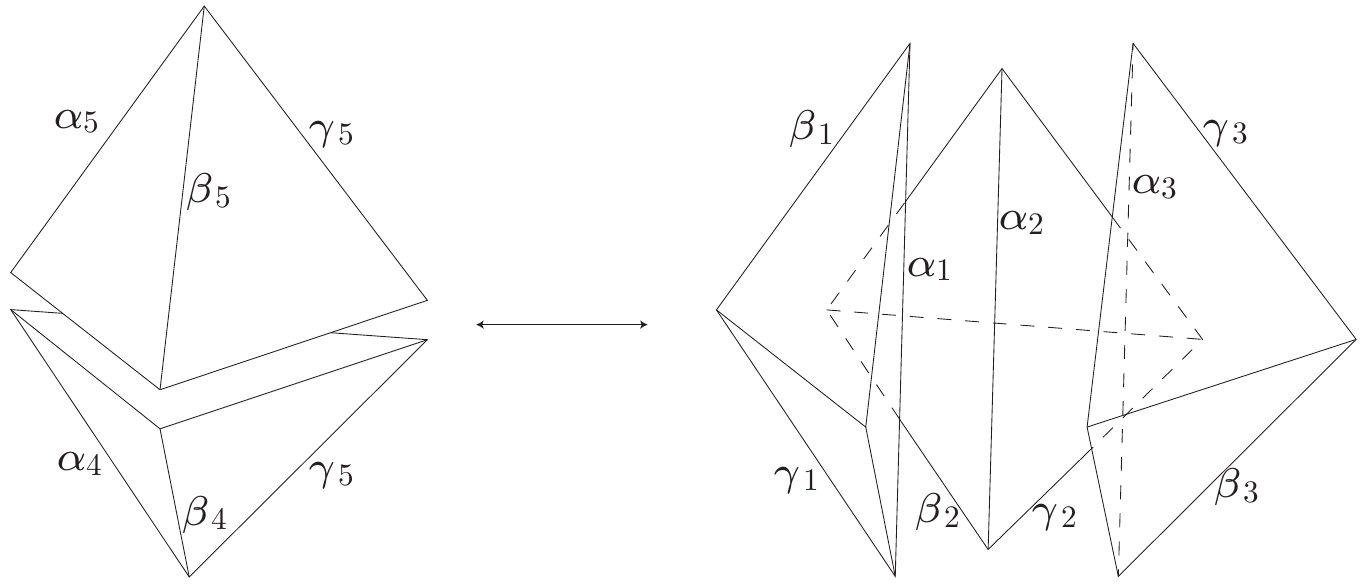}
		\caption{A $3$--$2$ Pachner move.}
		\label{fg:24}
	\end{figure}
	\begin{equation}\label{eq:Pchn}\begin{array}{ll}
	\alpha_4 = \beta_2 + \gamma_1 & \alpha_5 = \beta_1 + \gamma_2\\
	\beta_4 = \beta_1 + \gamma_3 & \beta_5 = \beta_3 + \gamma_1\\
	\gamma_4 = \beta_3 + \gamma_2 & \gamma_5 = \beta_2 + \gamma_3.
	\end{array}
	\end{equation}
	where $(\alpha_i,\beta_i,\gamma_i)$ are the dihedral angles of $t_i$.
	In this situation we say that $X_e$ is obtained from $X$ by a \emph{shaped $3-2$ Pachner move}.
\end{definition}

We remark that the linear system, together with $e$ being balanced, guaranties the positivity of the dihedral angles of $t_4$ and $t_5$ provided the positivity for $t_1$, $t_2$ and $t_3$ but it does not provide any guarantees on the converse, i.e.  the positivity of a shaped $2-3$ Pachner moves. However, two different solutions for the angles for $t_1$, $t_2$ and $t_3$ from the same starting angles for $t_4$ and $t_5$ are always gauge equivalent.\\
The system \eqref{eq:Pchn} define a map $P^e:\Shp(X) \longrightarrow \Shp(X_e)$, that extends to a map 
$$\tilde{P}^e:\GShp(X) \longrightarrow \GShp(X_e).$$ 
For a balanced edge $e$, the latter restricts to the map 
$$\tilde{P}_r : \tilde{\Omega}_{X,r}(e)^{-1} (2\pi)\longrightarrow \GShp_r(X_e),$$ 
and it can be noticed that $\tilde{P}_r( \tilde{ \Omega}_{X, r} (e)^{-1}(2\pi)  \cap \Shp_r(x)) \subset \Shp_r(Y)$.

We also say that a leveled shaped pseudo $3$-manifold $(X,l_X)$ is obtained from $(Y,l_Y)$ by a \emph{leveled shaped $3$-$2$ Pachner move} if, for some balanced $e\in\Delta_1(X)$, $Y=X_e$ as above and
$$ l_Y = l_X +\frac1{12\pi} \sum_{a\in(\phi^{3,1})^{-1} (e) }\sum_{b\in\Delta_3^1 (X)} \varepsilon_{p(a),p(b)} \alpha_X(b)\text{.}$$

\begin{definition}
	A (leveled)  shaped pseudo $3$-manifold $X$ is called a \emph{Pachner refinement} of a (leveled) shaped pseudo $3$-manifold $Y$ if there exists a finite sequence of (leveled) shaped pseudo $3$-manifolds
	\[
	X=X_1,X_2,\dots,X_n=Y
	\]
	such that for any $i\in\{1,\dots,n-1\}$, $X_{i+1}$ is obtained from $X_i$ by a (leveled) shaped $3 - 2$ Pachner move. Two (leveled) shaped pseudo $3$-manifolds $X$ and $Y$ are called \emph{equivalent} if there exist gauge equivalent (leveled) shaped pseudo 3-manifolds $X'$ and $Y'$ which are respective  Pachner refinements of $X$ and $Y$.
\end{definition}

For technical reasons, which we will discuss later, we will restrict the category of triangulated $2+1$ cobordisms, discussed so far, to a certain sub-categoroid, as discussed below.  This means that we will remove some morphisms as the following definition imposes. 

\begin{definition}[Admissibility]
	An oriented triangulated pseudo $3$-manifold is called \emph{admissible} if
	$$\Shp_r(X) \neq \emptyset\text{,}$$
	and
	$$
	\Hom_2(X-\Delta_0(X), \Integer) = 0.
	$$
\end{definition}

\begin{definition}
	Two (leveled) admissible shaped pseudo $3$-manifolds $X$ and $Y$ are said to be \emph{admissibly equivalent} if there exists a gauge equivalence $$h':X'\longrightarrow Y'$$ of (leveled) shaped $3$-manifolds $X'$ and $Y'$ which are respective Pachner refinements of $X$ and $Y$ such that $\Delta_1(X') = \Delta_1(X) \cup D_X$ and $ \Delta_1(Y') = \Delta_1(Y) \cup D_Y$ and the following holds
	$$\left[h(\Shp_r(X)\cap \tilde{\Omega}_{X',r}(D_X)^{-1}(2\pi))\right]\cap\left[  \tilde{\Omega}_{Y',r}(D_Y)^{-1}(2\pi)\right]\neq \emptyset\text{.}$$
\end{definition}

\begin{theorem}[\cite{AK1}]
	Suppose two (leveled) shaped pseudo $3$-manifolds $X$ and $Y$ are equivalent. Then there exist $D\subset \Delta_1(X)$ and $D' \subset \Delta_1(Y)$ and a bijection
	$$i : \Delta_1(X) - D \rightarrow  \Delta_1(Y) - D'$$
	and a Poisson isomorphism
	$$R :  \tilde \Omega_{X, r} (D)^{-1}(2\pi) \rightarrow \tilde \Omega_{Y, r} (D')^{-1}(2\pi),$$
	which is covered by  an affine $\Real$-bundle isomorphism from $\LGShp_r(X)|_{\tilde\Omega_{X,r}(D)^{-1} (2\pi)}$ to $\LGShp_r(Y)_{\tilde\Omega_{X,r}(D')^{-1} (2\pi)}$ and
	such that we get the following commutative diagram
	$$\begin{CD}
	\tilde \Omega_{X, r} (D)^{-1}(2\pi)  @>R>>  \tilde \Omega_{Y, r} (D')^{-1}(2\pi)\\
	@VV{\operatorname{proj}\ \circ \ \tilde \Omega_{X,r}}V @VV{\operatorname{proj}\ \circ \ \tilde \Omega_{Y,r}}V\\
	\Real^{\Delta_1(X)-D} @> i^* >>  \Real^{\Delta_1(Y)-D'}.
	\end{CD}
	$$
	Moreover, if $X$ and $Y$ are admissible and admissibly equivalent,  the isomorphism $R$ takes an open convex subset $U$ of $S_r(X)\cap \tilde \Omega_{X, r} (D)^{-1}(2\pi) $ onto an open convex subset $U'$ of  $S_r(Y)\cap \tilde \Omega_{Y, r} (D)^{-1}(2\pi) $.
\end{theorem}
We remark that in the previous notation $D= \Delta_1(X) \cap h^{-1} (D_Y)$ and $D' = \Delta_1(Y)\cap h(D_X)$.

For a tetrahedron $T=[v_0,v_1,v_2,v_3]$ in $\Real^3$ with ordered vertices $v_0,v_1,v_2,v_3$, we define its sign $$\sign(T)=\sign(\det(v_1-v_0,v_2-v_0,v_3-v_0)),$$
as well as the signs of its faces 
$$\sign(\partial_iT)=(-1)^{i}\sign(T)\text{,  for $ i\in\{0,\ldots,3\}$.}$$
For a pseudo 3-manifold $X$, the signs of faces of  the tetrahedra of $X$ induce a sign function on the faces of the boundary of $X$, $\sign_X: \Delta_2(\partial X)\to\{\pm1\}$, 
which permits us to split the components of the boundary of $X$ into two sets, $\partial X=\partial_+X\cup\partial_-X$, where  $\Delta_2(\partial_\pm X)=\sign_X^{-1}(\pm1)$.
Notice that $|\Delta_2(\partial_+ X)| = |\Delta_2(\partial_- X)|$.\newline
\begin{definition}[Cobordism Categoroid]
	The category $\Bordism$ is the category that has triangulated surfaces as objects,equivalence classes of (leveled) shaped pseudo $3$-manifolds $X$ as morphisms (so that $X\in\Homo_{\Bordism} (\partial_-X,\partial_+X)$) and the composition given by glueing along boundary components, through edge orientation preserving and face orientation reversing CW-homeomorphisms.\newline 
	The Categoroid $\Bordism_a$ is the subcategoroid of $\Bordism$ whose morphisms are restricted to be admissible equivalence classes of admissible (leveled) shaped pseudo $3$-manifolds. In particular composition is possible only if the gluing gives an (leveled) admissible pseudo $3$-manifold.
\end{definition} 

\begin{remark}\emph{Admissible Shaped Pseudo $3$-Manifolds in the real world.}\\
	Even though we will discuss the whole Andersen Kashaev construction of the Teichm\"uller TQFT functor in the general setting  
	of the  above defined cobordism categoroid, the the main parts of this construction, we want to put our hands on in this paper, are invariants of links. We  interpret Triangulated Pseudo $3$-manifolds $X$ as ideal triangulations of the (non closed) manifold $X\setminus \Delta_0(X)$. This interpretation is enlighten in Remark \ref{rm:idealTetra}. We shall ask ourself when a cusped $3$-manifold (cusped means non compact with finite volume here) admits a positive fully balanced shape structure. This requirement is weaker than asking for a full geometric structure on the manifold and in our language, this can be expressed by the fact that we did not required a precise gauge to be fixed. The problem of finding positive or generalized angle structures has been studied in \cite{LT}, where necessary and sufficient conditions for their existence are given. In the work \cite{HRS} it is proved, among other things, that a particular class of manifolds $M$ supporting positive shape structures are complements in $S^3$ of hyperbolic links. However the admissibility conditions kicks in here and further restrict us to just complements of hyperbolic knots. So, at the least, we know that the Andersen Kashaev construction will work on complements of hyperbolic knots, and that are the examples we will look a bit closer at below. Now we should clarify the equivalence relation in $\mathcal B_a$, in the context of knot complements. Combinatorially speaking, any two ideal triangulations of a knot complement are related by finite sequences of $3$--$2$ or $2$--$3$ Pachner moves. On the other hand it is not known (at least to the authors) if any such sequence can be realised as a sequence of shaped Pachner moves. For sure we know that $3$--$2$ shaped Pachner moves are well defined in the category $\mathcal B_a$ as we remarked when we defined them, and if a shaped Pachner $2$--$3$ move is possible in some particular case, then it is an equivalence in the category $\mathcal B_a$. So the knot invariants that we will define starting from $\mathcal B_a$ are not guaranteed  to be topological invariants. There is however another construction of the Andersen--Kashaev invariant \cite{AK2}, that avoid this problem with analytic continuation properties of the partition function. The equivalence of the two constructions is still conjectural though.\\
	In \cite{AK1}  another way to define knot invariants is suggested, by taking one vertex Hamiltonian triangulations of knots, that is, one vertex triangulations of $S^3$ (or a general manifold $M$) where the knot is represented by a unique edge with a degenerating shape structure, meaning that we take a limit on the shapes, sending all the weights to be balanced except the weight of the knot that is sent  to $0$. The partition function is actually divergent but a residue can be computed as an invariant. We will show this in a couple of examples in subsection \ref{sc:Comps}.
\end{remark}

\subsection{The target Categoroid $\D_N$}
W recall all the relevant things regarding tempered distributions and the space $\Schw(\Affine_N)$ in Appendix \ref{AP:distro}.  As always, here  $N$ is an odd positive integer and $\bpar\in\Complex$ is fixed to satisfy $\re(\bpar) >0$ and $\im\bpar(1-\abs{\bpar}) =0$.
\begin{definition}
	The categoroid $\D_N$ has as objects finite sets and for two finite sets $n,m$ the set of morphisms from $n$ to $m$ is
	$$
	\Homo_{\D_N}(n,m) = \Schw'(\Affine_N^{n\sqcup m})\simeq \Schw^{\p} (\Real^{n\sqcup m}) \otimes \Schw(\left(\IntN\right)^{n\sqcup m}).
	$$
\end{definition}

\begin{definition}
	For $\mathcal A\otimes A_N \in \Homo_{\D_N}(n,m)$ and $\mathcal B\otimes B_N\in \Homo_{\D_N}(m,l)$, such that $\mathcal A$ and $\mathcal B$  satisfy condition~(\ref{wftrans}) and  $\pi_{n,m}^*(\mathcal A)\pi_{m,l}^*(\mathcal B)$ continuously extends to $S(\Real^{n\sqcup m \sqcup l})_m$, we define
	$$
	(\mathcal A \otimes A_N)\circ( \mathcal B\otimes B_N) 
	= (\pi_{n,l})_*(\pi_{n,m}^*(\mathcal A)\pi_{m,l}^*(\mathcal B) )\otimes A_NB_N \in \Homo_{\D_N}(n,l).
	$$
\end{definition}
Where the product $A_NB_N$ is just the matrix product.\\
We will frequently use the following notation in what follows. For any $a\in\Affine_N$ , $a = (x,n)\in \Real\times\IntN$ we will consider the $\bpar$--dependent operator $$\varepsilon \equiv \varepsilon(\bpar) \colon \Affine_N \rightarrow \Affine_N$$ defined by
\begin{equation}\label{eq:epsilonDef}
\varepsilon (x,n) \equiv  \left\{ 
\begin{array}{l l}
(x,n) & \quad \text{if $\abs{\bpar} =1$,}\\
(x, -n) & \quad \text{if $\bpar\in\Real$ }\\
\end{array} \right.
\end{equation} 
For any $\mathcal A\in  \mathcal{L}( \Schw(\Real^{n}), \Schw'(\Real^{m}))$, we have a unique adjoint $\mathcal A^*\in \mathcal{L}( \Schw(\Real^{m}), \Schw'(\Real^{n}))$ defined by the formula
$$
\mathcal A^*(f)(g) = \overline{\bar{f}(\mathcal A(\bar{g}))}
$$ for all $f\in  \Schw(\Real^{m})$ and all $g\in \Schw(\Real^{n})$.

\begin{definition}[$\star_{\bpar}$ structure]\label{df:starFinite}
	Consider $\bpar\in\Complex$ fixed as above and $N\in\Integer_{>0}$ odd. Let $A_N\in\Homo(\Schw(\left(\IntN\right)^m), \Schw(\left(\IntN\right)^n)$. Recall the involution $\varepsilon$ on $\IntN$ form equation \eqref{eq:epsilonDef}. Define $A_N^{\star_{\bpar}}$ as
	\begin{align}
	\bra{j_1, \dots j_m}A_N^{\star_\bpar}\ket{p_1,\dots, p_n} = 
	\conj{\bra{\varepsilon p_1,\dots,\varepsilon p_n} A_N \ket{\varepsilon j_1, \dots\varepsilon j_m}   }
	\end{align}
\end{definition}
We can finally define the $*_\bpar$ operator as 
\begin{align}
\left(\mathcal A \otimes A_N\right)^{*_\bpar} = \mathcal A^*\otimes A_N^{\star_{\bpar}}
\end{align}

\subsection{Tetrahedral Partition Function}
Recall the operators from Section \ref{sc:ANreps}, $\opu_j$, $\opv_j$, $X_j$, $Y_j$, $\opp_j$, $\opq_j$ , $j =1,2$   acting on $\Hilb := \Schw(\Affine_N^2)$.
Define the \emph{Charged Tetrahedral Operator} as follows
\begin{definition}
	Let $a$, $b$, $c>0$ such that $a+b+c = \frac{1}{\sqrt N}$. Recall the Tetrahedral operator $\opT$ defined in \eqref{eq:TetraAdef}. Define the charged tetrahedral operator $\opT(a,c)$ as follows
	\begin{align}
	\opT(a,c) \equiv e^{-\pi i\frac{\cb^2}{\sqrt N} \left(2(a-c)+\frac{1}{\sqrt N}\right)/6}e^{2\pi i\cb( c\opq_2 - a\opq_1)} \opT_{12}
	e^{-2\pi i\cb( a\opp_2 + c\opq_2)}
	\end{align}
\end{definition}

\begin{lemma}\label{lm:CargedTetra}
	We have that \begin{equation}
	\opT(a,c) =  e^{-\pi i\frac{\cb^2}{\sqrt N} \left(2(a-c)+\frac{1}{\sqrt N}\right)/6} e^{\pi i \cb^2 a(a+c)}\opD_{12}\slashed{\psi}_{ a, c}(\opq_1+\opp_2-\opq_2,-e^{-\frac{\pi i}{N}}Y_1X_2\conj{Y_2})
	\end{equation}
	were $\psi_{a,c}(x,n)$ is the charged quantum dilogarithm from \eqref{eq:defCharDilog}
\end{lemma}

\emph{Extra Notation.}
Recall the notation for Fourier coefficients and Gaussian exponentials in $\Affine_N$. For $a = (x,n)$ and $a^\p = (y,m)$ in $\Affine_N$ we write
\begin{align*}
&\bkt{a,a^\p} \equiv e^{2\pi i xy}e^{-2\pi i nm/N}
&\bkt{a} \equiv e^{\pi i x^2} e^{-\pi i n(n+N)/N}
\end{align*}
For $a = (x,n)\in\Affine_N$, define $\delta(a)\equiv \delta(x) \delta(n)$ where $\delta(x)$ is Dirac's delta distribution while $\delta(n)$ is the Kronecker delta $\delta_{0,n}$ between $0$ and $n$ mod $N$. Define
\begin{equation}
\varphi_{a,c} (x,n)\equiv \psi_{a,c}(x, -n)\text{.}
\end{equation}
Denote, for $x, y\in\Real$ and $z\in\Affine_N$
\begin{align}
&\nu(x) \equiv e^{-\pi i\frac{\cb^2}{\sqrt N} \left(2x+\frac{1}{\sqrt N}\right)/6}  &\nu_{x,y}= \nu(x-y) e^{\pi i \cb^2 x(x+y)}
\end{align}
The equations from Lemma \ref{lm:dilogSym} can be upgraded to
\begin{align}
&\nu_{a,c}\tilde{\varphi}_{a,c} (z) = \nu_{c,b}\varphi_{c,b}(z) \bkt{z} e^{-\pi i N/12}\\
&\nu_{a,c}\conj{{\varphi}_{a,c}}(z) = \nu_{c,a}\varphi_{c,a}(-\varepsilon z) \bkt{z} e^{-\pi i N/6}\\
&\nu_{a,c}\conj{\tilde{\varphi}_{a,c}}(z) = \nu_{b,c} \varphi_{b,c}(-\varepsilon z) e^{-\pi i N/12}
\end{align}
Where $\varepsilon$ was defined in \eqref{eq:epsilonDef}.
From the  Charged Pentagon Equation \eqref{eq:charPenta} we get the following
\begin{proposition}[Charged Tetrahedral Pentagon equation]
	Let $a_j$, $c_j >0$ such that $\frac{1}{\sqrt{N}} - a_j - c_j >0$ for $j=0$, $1$, $2$, $3$ and $4$, which further satisfies the following relations 
	\begin{align}\label{eq:TpentaChargesCond}
	a_1 = a_0 +a_2\qquad a_3 = a_2 + a_4 \qquad c_1 = c_0 +a_4 \qquad c_3 = a_0 + c_4 \qquad c_2 = c_1 + c_3\text{.}
	\end{align}
	Then we have that
	\begin{align}\label{eq:ChTetraPenta}
	\opT_{12}(a_4,c_4)\opT_{1,3} (a_2, c_2 ) \opT_{23} (a_0, c_0) = \mu \opT_{23}(a_1, c_1)\opT_{12}(a_3,c_3)
	\end{align}
	where
	$$\mu = \exp  \pi i \frac{\cb^2}{6\sqrt N} \left(2(c_0 +a_2 + c_4) - \frac{1}{\sqrt N} \right)$$
\end{proposition}

We have an integral kernel description for the charged tetrahedral operator. We use the Bra-Ket notation to denote integral kernels, see Appendix \ref{ap:braket}.
\begin{proposition}\label{pr:TKernel}
	Let $\conj{\opT}(a,c) \equiv (\opT(a,c))^{-1}$.
	\begin{align*}
	&\bra{a_0,\, a_2}\, \opT_{12}(a,c) \,\ket{a_1, \, a_3} \\
	&\qquad=\nu(a-c) e^{\pi i \cb^2 a(a+c)}\bkt{a_3 - a_2,\, a_0} \conj{\bkt{a_3 - a_2}}\delta(a_0+a_2-a_1) \tilde{\varphi}_{a,c}(a_3 -a_2)\\
	&\bra{a_0, \, a_2 } \conj{\opT}(a,c) \ket{a_1, \, a_3} \\
	&\qquad=\nu(b\!-\!c) e^{\pi i \cb^{2} b(b+c)}e^{-\pi i N/12} \bkt{a_3\! -\! a_2, a_1} \bkt{a_3\! -\! a_2}\delta(a_1\! +\! a_3\! -\! a_0)\varphi_{b,c} (a_3\! -\! a_2).
	\end{align*}
\end{proposition}

The appearance of $\varepsilon$ is due to the non-unitarity of the theory for $\bpar >0$ and $N>1$.

Let $\opA$ and $\opB$ two operators on $\SqInt(\Affine_N)$ defined as bra-ket distributions by
\begin{align}
&\bkt{a_1, a_2| \opA} = \delta(a_1 + a_2)  \bkt{a_1}e^{\pi i N/12}   & \bkt{a_1, a_2| \opB} =   \bkt{a_1 -a_2}\\
&\bkt{\conj{\opA}| a_1, a_2} = \conj{\bkt{\varepsilon a_1,\varepsilon a_2| \opA}}     &\bkt{\conj{\opB}| a_1, a_2} =  \conj{\bkt{\varepsilon a_1, \varepsilon a_2| \opB}}.
\end{align}
\begin{lemma}[Fundamental Lemma]\label{lm:FundLem}
	We have the following three relations
	\begin{align}
	&\int_{\Affine_N^2}\bkt{\conj{\opA}| v,s}\bra{x,s}\opT(a,c)\ket{u, t}  \bkt{t,y| \opA}\diff s \diff t = 
	\bra{x,y}\conj\opT(a,b)\bkt{u, v}\\
	&\int_{\Affine_N^2}\bkt{\conj{\opA}| u,s}\bra{s,x}\opT(a,c)\ket{v, t}  \bkt{t,y| \opB}\diff s \diff t = 
	\bra{x,y}\conj\opT(b,c)\bkt{u, v}\\
	&\int_{\Affine_N^2}\bkt{\conj{\opB}| u,s}\bra{s,y}\opT(a,c)\ket{t, v}  \bkt{t,x| \opB}\diff s \diff t = 
	\bra{x,y}\conj\opT(a,b)\bkt{u, v}.
	\end{align} 
\end{lemma}

\subsubsection{TQFT Rules, Tetrahedral Symmetries and Gauge Invariance}
We consider oriented surfaces with  cellular structure such that all $2$-cells are either bigons or triangles. Not all the edge orientations will be admitted. We forbid cyclically oriented triangles.
For the bigons, we consider only the \emph{essential} ones, the others being contractible to an edge. These essential bigons are precisely the ones with cellular structure isomorphic to the unit disk with vertices $\pm 1\in\Complex$ and edges  $\{e_1 = e^{\pi i t};\, e_2 = -e^{\pi it} \text{, for } t\in[0,1]\}$ or $\{e_1 =-e^{-\pi i t};\, e_2 = e^{-\pi it}\text{, for } t\in [0,1]\}$.
Given such an ideally triangulated surface $\Sigma$ we will associate a copy of $\Complex$ to any bigon and a copy of $\Schw'(\Affine_N)$ to any triangle. Globally we  associate  to the surface the space $\Schw'(\Affine_N^{\Delta_2(\Sigma)})$.
To a shaped tetrahedron $T$ with ordered vertices $\{v_0,v_1,v_2,v_3\}$ we associate the partition function $Z^{(N)}_{\bpar} (T)$ through the Nuclear Theorem \eqref{eq:NuclearBraKet} as a ket distribution
\begin{equation}\label{eq:PFtetra}
\bkt{ x|\widetilde{Z_{\bpar}^{(N)} (T)}}=
\left\{
\begin{array}{cc}
\bra{ a_0,a_2}\opT(c(v_0v_1),c(v_0v_3))\ket{ a_1,a_3}&\mathrm{if}\ \sign(T)=1;\\
\bra{ a_1,a_3}\conj\opT(c(v_0v_1),c(v_0v_3))\ket{ a_0,a_2}&\mathrm{if}\ \sign(T)=-1.
\end{array}
\right.
\end{equation}
where
$$\Affine_N\,\ni \,a_i := a(\partial_i T)\text{, \hskip 1em} i\in\{0,1,2,3\}$$
and
$$c:= \frac 1{\pi\sqrt N} \alpha_T :\Delta_1(T) \rightarrow \Real_{>0}\text{.}$$

Having allowed bigons in triangulations of surfaces, we must also allow cones over such as cobordisms. From the $2$ classes of bigons described above we have $4$ isotopy classes of cellular structures of cones over them, described in the following as embedded in $\Real^3\simeq \Complex\times\Real$. The bigon is identified with the unit disc embedded in $\Complex$.  The apex of the cone will be the point $(0,1)\in \Complex\times\Real$. The $1$-cells will be either
\[
\{e^1_{0\pm}(t)=(\pm e^{i\pi t},0),\ e^1_{1\pm}(t)=(\pm(1-t),t)\}
\]
or
\[
\{e^1_{0\pm}(t)=(\mp e^{-i\pi t},0),\ e^1_{1\pm}(t)=(\pm(1-t),t)\}
\]
or
\[
\{e^1_{0\pm}(t)=(\pm e^{i\pi t},0),\ e^1_{1\pm}(t)=(\pm t,1-t)\}
\]
or
\[
\{e^1_{0\pm}(t)=(\mp e^{-i\pi t},0),\ e^1_{1\pm}(t)=(\pm t,1-t)\}.
\]
We name these types of cones $A_+$, $A_-$, $B_+$  and $B_-$ respectively. We need TQFT rules for the gluing of these cones. We just need to consider their gluing to a tetrahedra.
We assign a partition function to the cones as follows
\begin{equation}\label{eq:SymmPartFunc}
\bkt{a_1,a_2| \widetilde{Z_{\bpar}^{(N)}(A_{\pm})}} = \delta(a_1 + a_2)  \bkt{a_1}^{\pm 1}e^{\pm\pi i N/12} \text{, \hskip 1em} \bkt{a_1,a_2| \widetilde{Z_{\bpar}^{(N)}(B_{\pm})}} =   \bkt{a_1 -a_2}^{\pm 1}\text{.}
\end{equation}
Tetrahedral symmetries are generated by permutation of the ordered vertices. Indeed the group of tetrahedral symmetries is identified with the symmetric group $\mathbb S_4$ and is generated by three transpositions. The three equations of the Fundamental Lemma \ref{lm:FundLem} gain an interpretation as glueing of cones on the faces of a tetrahedron through  definitions \eqref{eq:SymmPartFunc}. These three glueing generates all the symmetries of a tetrahedron, and through this interpretation, the Fundamental Lemma assure that the partition function $Z_{\bpar}^{(N)}$ satisfies all the tetrahedral symmetries. For more details on tetrahedral symmetries and the cone's partition function see \cite{AK1, GKT}.

We can now formulate the main Theorem for the Teichm\"uller TQFT. This theorem was proved by Andersen and Kashaev for the case $N=1$ in \cite{AK1}. The statement that we have here is for every $N$ odd, and it is strictly speaking not present as such in the literature.
\begin{theorem}[Level $N$ Teichm\"uller TQFT, Andersen and Kashaev]\label{th:TeichTQFT}
	For any $\bpar \in\Complex^*$ such that  $\im \bpar(\abs{\bpar} - 1) =0$ and $\re\bpar >0$, and for any $N\in\Integer_{>0}$ odd there exists a unique $*_{\bpar}$-functor
	$F_{\bpar}^{(N)}\colon \Bordism_a \rightarrow \D_N$ such that
	\(
	F_\bpar^{(N)}(A)=\Delta_2(A),\  \forall A\in\operatorname{Ob}\Bordism_a,
	\)
	and for any admissible leveled shaped  pseudo 3-manifold $(X,l_X)$, the associated morphism in $\D_N$ takes the form
	\begin{equation}\
	F_\bpar^{(N)}(X,l_X)=Z_{\bpar}^{(N)}(X)e^{-\pi i\frac{l_X\cb^2}{N}}\in\Schw'\left(\Affine_N^{\Delta_2(\partial X)}\right),
	\end{equation}
	where $Z_{\bpar}^{(N)}$ is defined in \eqref{eq:PFtetra} for a tetrahedron.
\end{theorem}

Here $*_{\bpar}$-functor means that $F_\bpar^{(N)}(X^*)=F_\bpar^{(N)}(X)^{*_{\bpar}}$, where $X^*$ is the oppositely oriented pseudo $3$-manifold to $X$.

The discussion so far proves the theorem except for the gauge invariance and the convergence of the partition functions under glueings. We will not discuss the convergence here because it follows directly from the convergence in the case level $N=1$, which was addressed in \cite{AK1}. We just remark that the hypothesis of admissibility is used to prove the convergence of the partition function.\\
For the gauge invariance consider the suspension of an $n$-gone $SP_n$ naturally triangulated into $n$ tetrahedra sharing the only internal edge $e$. Every gauge transformation can be decomposed in a sequence of gauge transformations involving only one edge $e$, and every such gauge transformation can be understood in the example of the suspension.
Suppose all the tetrahedra to be positive, and having vertex order such that the last two vertices are the endpoints of the internal common edge. After enumerating the tetrahedra in cyclic order, let $a_i$, $c_i$ be the two shape parameter of $T_i$, $i=0,\dots,n$, and $\mathrm{a} = (a_0,\dots , a_n)$, $\mathrm{c} = (c_0,\dots,c_n)$. Notice that $\sqrt N \pi a_i$ is the dihedral angle corresponding to the edge $e$. So a gauge tranformation corresponding to $e$ will affect the partition function of $SP_N$ 
$$Z_{\bpar}^{(N)} (SP_N)(\mathrm{a},\mathrm{c}) := \operatorname{Tr}_0(\opT_{01}(a_1,c_1)\opT_{02} (a_2,c_2)\cdots\opT_{0n}(a_n,c_n))$$ 
by shifting $\mathrm{c}$ by an amount $\mathrm{\lambda} = (\lambda,\dots,\lambda)$ say. One can show from the definitions and the discussion above, 
that 
$$\opT(a,c+\lambda) =e^{-2\pi i \cb\lambda\opp_1} \opT(a,c)  e^{2\pi i \cb\lambda\opp_1} e^{\pi i \cb^2 \left(\frac{1}{\sqrt N}- 6a\right)\lambda/3}$$
which, after tracing, leads to the following 
\begin{proposition}\cite{AK1}
	$$Z_{\bpar}^{(N)} (SP_N)(\mathrm{a},\mathrm{c} + \mathrm{\lambda}) = Z_{\bpar}^{(N)} (SP_N)(\mathrm{a},\mathrm{c}) e^{\pi i\cb^2\left(\frac{n}{\sqrt N}-6Q_e\right)\lambda/3}$$
	where 
	$$Q_e = a_1+a_2 + \dots a_n$$
\end{proposition}

\subsection{Knot Invariants: Computations and Conjectures}\label{sc:Comps}

In this secttion we update the examples computed in \cite{AK1} to the level $N\geq 1$ setting. Similar results were obtained in \cite{D3D}.\\
\emph{Notation.}  
In the examples we are going to use the following notation for quantum dilogarithms
\begin{align}
\varphi_{\bpar} (x,n) \equiv \bDilog(x, -n).
\end{align}
Moreover we will often abuse notation in favor of readability in the following ways. For $z=(x,n)\in\Affine_N$ we will sometimes write
$e^{2\pi i \cb z \alpha}$ in place of $e^{2\pi i \cb x\alpha}$. Moreover sums of the form $z+\cb a$ will always mean $(x+\cb a, n)$. 

In the following examples we encode an oriented triangulated pseudo 3-manifold $X$ into a diagram where a tetrahedron $T$ is represented by an element
\begin{center}
	\begin{tikzpicture}[scale=.3]
	\draw[very thick] (0,0)--(3,0);
	\draw (0,0)--(0,1);\draw (1,0)--(1,1);\draw (2,0)--(2,1);\draw (3,0)--(3,1);
	\end{tikzpicture}
\end{center}
where the vertical segments, ordered from left to right, correspond to the faces $\partial_0T,\partial_1T,\partial_2T,\partial_3T$ respectively. When we glue tetrahedron along faces, we illustrate this by joining the corresponding vertical segments.
\subsubsection{Figure--Eight  Knot $4_1$}
Let $X$ be represented by the diagram
\begin{equation}\label{E:8graph}
\begin{tikzpicture}[baseline=5pt,scale=.3]
\draw[very thick] (0,0)--(3,0);\draw[very thick] (0,1)--(3,1);
\draw(0,0)--(1,1);
\draw (1,0)--(0,1);
\draw (2,0)--(3,1);
\draw (3,0)--(2,1);
\end{tikzpicture}
\end{equation}
Choosing an orientation,
it consists of one positive tetrahedron $T_+$ and one negative tetrahedron $T_-$ with four identifications
\[
\partial_{2i+j}T_+\simeq\partial_{2-2i+j}T_-,\quad i,j\in\{0,1\}\text{.}
\]
Combinatorially, we have $\Delta_0(X)=\{*\}$, $\Delta_1(X)=\{e_0,e_1\}$, $\Delta_2(X)=\{f_0,f_1,f_2,f_3\}$, and $\Delta_3(X)=\{T_+,T_-\}$ with the boundary maps
\[
f_{2i+j}=\partial_{2i+j}T_+=\partial_{2-2i+j}T_-,\quad i,j\in\{0,1\},
\]
\[
\partial_if_j=\left\{\begin{array}{cl}
e_0,&\mathrm{if}\ j-i\in\{0,1\};\\
e_1,&\mathrm{otherwise},
\end{array}
\right.
\]
\[
\partial_i e_j=*,\quad i,j\in\{0,1\}.
\]
The topological space $X\setminus\{*\}$ is homeomorphic to the complement of the figure--eight knot, and indeed $X\setminus \{*\}$ is an ideal triangulation of such a cuspidal manifold. The set $\Delta_{3,1}(X)$ consists of the elements $(T_\pm, e_{j,k})$ for $0\le j<k\le 3$. We fix a shape structure
\[
\alpha_X\colon \Delta_{3,1}(X)\to \Real_{>0}
\]
by the formulae
\[
\alpha_X(T_\pm, e_{0,1})=\pi\sqrt N a_\pm,\quad \alpha_X(T_\pm, e_{0,2})=\pi \sqrt N b_\pm,\quad\alpha_X(T_\pm, e_{0,3})=\pi \sqrt N c_\pm,
\]
where $a_\pm+b_\pm+c_\pm=\frac{1}{\sqrt N}$. The weight function
\[
\omega_X\colon \Delta_1(X)\to\Real_{>0}
\]
takes the values
\[
\omega_X(e_0)=\sqrt N\pi(2a_++c_++2b_-+c_-)=:2\pi w,\quad \omega_X(e_1)=2\pi(2-w).
\]
As the figure--eight knot is hyperbolic, the completely balanced case $w=1$ is accessible directly. We can state the balancing condition $w = 1$ as 
\begin{equation}
2b_+ + c_+ = 2b_- + c_-\text{.}
\end{equation}

The kernel representations for the operators $\opT(a_+, c_+)$ and $\opT(a_-, c_-)$ are as follows. Let $z_j\in\Affine_N$ , $j =0,\,1,\,2,\,3$,
\begin{align}
&\bra{z_0, z_2}\opT(a_+,c_+) \ket{z_1, z_3} \\
\nonumber
&\qquad \qquad = \nu_{a_+,c_+}\bkt{z_3-z_2, z_0} \conj{\bkt{z_3 - z_2} }  \delta(z_0 +z_2 - z_1) \tilde{\varphi}_{a_+, c_+}(z_3 - z_2)\\
&\bra{z_3, z_1}\conj{\opT}(a_+,c_-) \ket{z_2, z_0} = \conj{\bra{\varepsilon z_2,\varepsilon z_0}{\opT}(a_+,c_-) \ket{\varepsilon z_3,\varepsilon z_1}}\\
\nonumber
&\qquad \qquad = \conj{\nu_{a_-,c_-}}\bkt{z_0-z_1, z_2} {\bkt{z_1 - z_0} }  \delta(z_0 +z_2 - z_3) \conj{\tilde{\varphi}_{a_-, c_-}(\varepsilon z_1 - \varepsilon z_0)}
\end{align}

The Andersen--Kashaev invariant at level $N$ for the complement of the figure--eight knot is then
\begin{align*}
Z^{(N)}_{\bpar} (X) &= \int_{\Affine_N^4} \bra{z_0, z_2}\opT(a_+,c_+) \ket{z_1, z_3}\bra{z_3, z_1}\conj{\opT}(a_+,c_-) \ket{z_2, z_0}  \diff z_0 \diff z_1 \diff z_2 \diff z_3\\
&= \int_{\Affine_N^4}\nu_{c_+,b_+} \conj{\nu_{c_-,b_-}}{\varphi}_{c_+, b_+}(z_3 - z_2)\conj{{\varphi}_{c_-, b_-}(\varepsilon z_1 - \varepsilon z_0)}   \delta(z_0 +z_2 - z_1) \times\\
&\qquad\times \delta(z_0 +z_2 - z_3) \bkt{z_3-z_2, z_0}\bkt{z_0-z_1, z_2}\diff z_0 \diff z_1 \diff z_2 \diff z_3\\
&= \int_{\Affine_N^3}\nu_{c_+,b_+} \conj{\nu_{c_-,b_-}}{\varphi}_{c_+, b_+}(z_1 - z_2)\conj{{\varphi}_{c_-, b_-}(\varepsilon z_1 - \varepsilon z_0)}   \delta(z_0 +z_2 - z_1) \times\\
&\qquad\times \bkt{z_1-z_2, z_0}\bkt{z_0-z_1, z_2}\diff z_0 \diff z_1 \diff z_2 \\
&= \int_{\Affine_N^2}\nu_{c_+,b_+} \conj{\nu_{c_-,b_-}}{\varphi}_{c_+, b_+}(z_0)\conj{{\varphi}_{c_-, b_-}(\varepsilon z_2)}    \bkt{z_0, z_0}\bkt{-z_2, z_2}\diff z_0  \diff z_2 \\
&= \int_{\Affine_N}\nu_{c_+,b_+} {\varphi}_{c_+, b_+}(z_0)\bkt{z_0, z_0}\diff z_0
\conj{\int_{\Affine_N} \nu_{c_-,b_-} {{\varphi}_{c_-, b_-}(\varepsilon z_2)}\bkt{z_2, z_2}  \diff z_2 }\\
&= \sigma_{c_+,b_+}\conj{\sigma_{c_-,b_-}}
\end{align*}
We can compute 
\begin{align*}
\sigma_{c_\pm, b_\pm} &= \nu_{c_\pm,b_{\pm}}\int_{\Affine_N}\frac{e^{-2\pi i \cb zc_{\pm}}}{\varphi_{\bpar} (z - \cb(b_\pm +c_\pm))}\bkt{z}^2\diff z \\
&= \nu^{\p}_{c_\pm, b_\pm} \int_{\Affine_N + di}\frac{e^{4\pi i \cb z(2b_{\pm} +c_{\pm})}}{\varphi_{\bpar} (z )}\bkt{z}^2\diff z \\
\end{align*}
where
\begin{align}
\nu^{\p}_{c_\pm, b_\pm}  = \nu_{c_\pm,b_{\pm}} e^{4\pi i \cb^2(c_{\pm}b_{\pm} - b_{\pm}^2)}
\end{align}
and the domain of integration $\Affine_N + di = \left(\Real+di)\times\IntN\right)$. Note we have shifted the real integral to a contour integral in the complex plane, and $d\in\Real$ is such that the integral converges absolutely. We sometimes omit the contour shift in the computations but we state it in the results.
Defining
\begin{align*}
\lambda \equiv 2b_+ +c_+ = 2b_- + c_- 
\end{align*}
we have

\begin{align*}
\begin{split}
Z^{(N)}_{\bpar} (X) &=\nu^{\p}_{c_+, b_+}\conj{\nu^{\p}_{c_-, b_-}} \int_{\Affine_N^2}\frac{e^{4\pi i\cb \lambda(z_0+z_2)}}{\varphi_{\bpar}(z_0) \conj{\varphi_{\bpar}(\varepsilon z_2)}} \bkt{z_0}^2 \conj{\bkt{z_2}^2}\diff z_0 \diff z_2\\
&= \nu^{\p}_{c_+, b_+}\conj{\nu^{\p}_{c_-, b_-}} \int_{\Affine_N^2}\frac{{\varphi_{\bpar}( z_2)}}{\varphi_{\bpar}(z_0)}  e^{4\pi i\cb \lambda(z_0+z_2)}\bkt{z_0}^2 \conj{\bkt{z_2}^2}\diff z_0 \diff z_2\\
&= \nu^{\p}_{c_+, b_+}\conj{\nu^{\p}_{c_-, b_-}} \int_{\Affine_N^2}\frac{{\varphi_{\bpar}( z_2 - z_0)}}{\varphi_{\bpar}(z_0)}  e^{4\pi i \cb\lambda z_2}\bkt{z_0, z_2}^2 \conj{\bkt{z_2}^2}\diff z_0 \diff z_2
\end{split}
\end{align*}
that has the structure
\begin{align}
&Z^{(N)}_{\bpar} (X) = e^{i\phi}  \int_{\Affine_N+i0} \chi^{(N)}_{4_1}(x, \lambda)   \diff x\text{,}\\
&\chi^{(N)}_{4_1}(x, \lambda) = \chi^{(N)}_{4_1}(x) e^{4\pi i\cb \lambda x}\text{,} \qquad\qquad 
\chi^{(N)}_{4_1}(x) =\int_{\Affine_N-i0}\frac{{\varphi_{\bpar}(x - y)}}{\varphi_{\bpar}(y)} \bkt{x, y}^2 \conj{\bkt{x}^2}\diff y
\end{align}
where $\phi$ is some constant quadratic combination of dihedral angles.

\subsubsection{The Complement of the Knot $5_2$}
Let $X$ be the closed S.O.T.P. $3$-manifold represented by the diagram
\[
\begin{tikzpicture}[scale=.3]
\draw[very thick] (0,0)--(3,0);
\draw[very thick] (6,0)--(9,0);
\draw[very thick] (3,3)--(6,3);
\draw(3,0)..controls (3,1) and (6,1)..(6,0);
\draw(0,0)..controls (0,2) and (3,1)..(3,3);
\draw(1,0)..controls (1,2) and (4,1)..(4,3);
\draw(2,0)..controls (2,2.5) and (9,2.5)..(9,0);
\draw(5,3)..controls (5,1) and (8,2)..(8,0);
\draw(6,3)..controls (6,2) and (7,2)..(7,0);
\end{tikzpicture}
\]
This triangulation has only one vertex $*$ and $X\setminus\{*\}$ is topologically the complement of the knot $5_2$. 
We denote $T_1,T_2,T_3$ the left, right, and top tetrahedra respectively. We choose the orientation so that  all of them are positive.
Balancing all the edges correspond to require the following equations to be true
\begin{align}\label{eq:balan52}
2a_3=a_1+c_2,\quad b_3=c_1+b_2.
\end{align}
The three integral kernels reads
\begin{align*}
&\bra{z,w} \opT(a_1, c_1)\ket{u,x} =\\
&\qquad\qquad =\nu_{a_1,c_1}\bkt{x-w, z}\conj{\bkt{x-w}} \delta(z+w-u)\tilde{\varphi}_{a_1,c_1}(x-w)\\
&\bra{x,v} \opT(a_2, c_2)\ket{y,w} =\\
&\qquad\qquad =\nu_{a_2,c_2}\bkt{w - v, x}\conj{\bkt{w-v}} \delta(x+v -y)\tilde{\varphi}_{a_2,c_2}(w - v)\\
&\bra{y,u} \opT(a_3, c_3)\ket{v,z} =\\
&\qquad\qquad =\nu_{a_3,c_3}\bkt{z-u, y}\conj{\bkt{z-u}} \delta(y+u -v)\tilde{\varphi}_{a_3,c_3}(z-u)\\
\end{align*}
Carrying out the computations, defining $\lambda =  -c_1 + b_2 -c_2 + a_3$, one gets that
\begin{align}
&Z^{(N)}_{\bpar} (X) = \int_{\Affine_N+i0}\chi_{5_2}^{(N)}(x,\lambda)\diff x\text{,}  \qquad \qquad
\chi_{5_2}^{(N)}(x,\lambda) = \chi_{5_2}^{(N)}(x)e^{2\pi i \cb\lambda x}\\
&\chi_{5_2}^{(N)}(x) = \int_{\Affine_N-i0}\frac{\conj{\bkt{x }} \bkt{ z} }
{\varphi_{\bpar}(z+x)\varphi_{\bpar}(z)\varphi_{\bpar}(z-x)}\diff z
\end{align}

\subsubsection{H-Triangulations}
In this section we will look at one vertex H-triangulations of knots.

Let $X$ be an H--Triangulation for the figure--eight knot, i.e. let $X$ be given by the diagram
\[
\begin{tikzpicture}[scale=.3]
\draw[very thick] (0,0)--(0,3);
\draw[very thick] (1,3/2)--(4,3/2);
\draw[very thick] (5,0)--(5,3);
\draw(1,3/2)..controls (1,2) and (2,2)..(2,3/2);
\draw(3,3/2)..controls (3,2) and (3/2,5/2)..(0,2);
\draw(4,3/2)..controls (4,2) and (4.5,3)..(5,3);
\draw(0,3)..controls (1/2,3) and (4.5,2)..(5,2);
\draw(0,1)..controls (1/2,1) and (4.5,0)..(5,0);
\draw(0,0)..controls (1/2,0) and (4.5,1)..(5,1);
\end{tikzpicture}
\]
where the figure-eight knot is represented by the edge of the central tetrahedron connecting the maximal and the next to maximal vertices.
If we choosing central tetrahedron ($T_0$) to be positve, the left tetrahedron ($T_+$) will be positive and the right one ($T_-$) negative.The shape structure, in the limit $a_0\rightarrow 0$ satisfies $2b_++c_+=2b_-+c_-=:\lambda$
The partition function satisfies the following limit formula 

\begin{align}
\lim_{a_0 \rightarrow 0}\varphi_{\bpar}(\cb a_0 -\cb/\sqrt N)Z_\bpar^{(N)}(X) =\frac{e^{-\pi i N/12}}{\nu(c_0)} \chi_{4_1}^{(N)}(0)
\end{align}

Similarly let $X$ be represented by the diagram
\[
\begin{tikzpicture}[scale=.3]
\draw[very thick] (0,0)--(3,0);
\draw[very thick] (3,1)--(6,1);
\draw[very thick] (6,0)--(9,0);
\draw[very thick] (3,3)--(6,3);
\draw(5,1)..controls (5,1/2) and (6,1/2)..(6,1);
\draw (3,0)--(3,1);
\draw(4,1)..controls (4,0) and (6,.5)..(6,0);
\draw(0,0)..controls (0,2) and (3,1)..(3,3);
\draw(1,0)..controls (1,2) and (4,1)..(4,3);
\draw(2,0)..controls (2,2.5) and (9,2.5)..(9,0);
\draw(5,3)..controls (5,1) and (8,2)..(8,0);
\draw(6,3)..controls (6,2) and (7,2)..(7,0);
\end{tikzpicture}
\]
that is, the H--triangulation for the $5_2$ knot.
We denote $T_0,T_1,T_2,T_3$ the central, left, right, and top tetrahedra respectively and we choose the orientation so that the central tetrahedron $T_0$ is negative then all other tetrahedra are positive.
The edge representing the knot $5_2$ connects the last two edges of  $T_0$, so that the weight on the knot is given by $2\pi a_0$.  In the limit $a_0\to0$,  all edges, except for the knot, become balanced under the conditions
\[
a_1=c_2=a_3,\quad b_3=c_1+b_2,
\]
which in particular imply \eqref{eq:balan52}. 
The partition function has the following expression
\begin{align}
Z_{\bpar}^{(N)}(X) 
&= \Theta \frac{e^{-\pi i N/12}}{\varphi_{\bpar}(\cb a- \cb \sqrt N)}\chi_{5_2}^{(N)}(\cb(a_1-a_3))
\end{align}
For some constant phase factor $\Theta$.

\subsubsection{Asymptotic's of $\chi_{4_1}^{(N)}(0)$} \label{sc:AsyKnot}

In this section we want to study the asymptotic behavior of the invariant of the figure--$8$ knot 
\begin{align*}
\begin{split}
\chi_{4_1}^{(N)}(0) &= \int_{\Affine_N} \thDilog(-x, - k) \conj{\thDilog(x,-k)}\diff (x,k)\\
&= \frac{1}{2\pi \bpar\sqrt{N} }\sum_{k\in\IntN} \int_{\Real - id}\thDilog\left(\frac{-x}{2\pi \bpar}, - k\right) \conj{\thDilog\left(\frac{x}{2\pi \bpar},-k\right)}\diff x 
\end{split}
\end{align*}
when $\bpar \rightarrow 0$. The analysis uses techniques similar to the one presented in \cite{AK1} for $N=1$, however higher level gives new informations that we will show here.\\
The integration in the complex plane is a contour integral where $d>0$ so that the integral is absolutely convergent.
By means of the asymptotic formula for the quantum dilogarithm \eqref{eq:AsymFaddeev} we have that
\begin{align*}
\chi_{4_1}^{(N)}(0) =&\frac{1}{2\pi \bpar\sqrt{N} }\sum_{k\in\IntN} \int_{\Real - id}\operatorname{Exp}\left[\frac{\Dilog(-e^{-\sqrt Nx}) - \Dilog(-e^{\sqrt Nx})}{2\pi i\bpar^2 N} \right]\\
&\qquad\times  \phi_{-x}(k) \conj{\phi_x(k)}(1+\mathcal O(\bpar^2))\diff x
\end{align*}
We want to apply the steepest descent method to this integral to get an asymptotic formula for $\bpar \rightarrow 0$.
First we show the computation for the exact integral,
\begin{equation}\label{simpleIntegral}
\frac{1}{2\pi \bpar\sqrt{N} }\sum_{k\in\IntN} \int_{\Real - id}\operatorname{Exp}\left[\frac{\Dilog(-e^{-\sqrt Nx}) - \Dilog(-e^{\sqrt Nx})}{2\pi i\bpar^2 N} \right] \phi_{-x}(k) \conj{\phi_x(k)}\diff x
\end{equation}
and then we will argue that the former one can be approximated by the  latter when $\bpar \rightarrow 0$.\\
Let $h(x) := \Dilog(-e^{-\sqrt Nx}) - \Dilog(-e^{\sqrt Nx})$. Its critical points are solutions to 
$$\begin{cases}h^{\prime}(x) = 0\\
h^{\prime\prime}(x) \neq 0
\end{cases}$$
which are $\mathcal S = \left\{ \pm\frac23\frac{\pi i }{\sqrt N} +\frac{2\pi i k}{\sqrt N} : k\in\Integer\right\}$. 
We compute the value of $\im h$ at its critical points to be
\begin{equation}
\im h\left(\pm\frac23\frac{\pi i }{\sqrt N} + \frac{2\pi i k}{\sqrt N}\right) = \pm4\Lambda(\frac\pi 6)
\end{equation}
where $\Lambda$ is the Lobachevsky's function
\begin{equation}
\Lambda(\alpha) = -\int_0^{\alpha } \log\abs{2\sin\varphi}\diff\varphi
\end{equation} 
and we refer the reader to \cite{KirillovDilogIdentities} for the expressions that relate Lobachevsky's function to the classical dilogarithm.

We only remark  that $4\Lambda(\frac{\pi}{6}) =  \operatorname{Vol}(4_1)$, where by $\operatorname{Vol}(4_1)$ we mean the hyperbolic volume of knot complement $S^3\setminus (4_1)$.\\
Fix $\Complex\ni x_0= -\frac23\frac{\pi i}{\sqrt N}$, which is accessible from the original contour without passing through other critical points, and consider  the contour 
$$\mathcal C = \left\{ z\in\Complex : \re(h(z)) = \re(h(x_0))\text{, } \im(h(z)) \leq \im(h(x_0)) \right\} $$
which is asymptotic to $\re(z) +\im(z) = 0$ for $\re(z) \rightarrow \infty$ and 
to $\re(z) -\im(z) = 0$ for $\re(z) \rightarrow -\infty$. Moreover 
\begin{equation} 
\lim_{\re(z)\rightarrow \pm\infty} \im( h(z)) = \lim_{\re(z)\rightarrow \pm\infty} \pm\re(z)\im(z) = -\infty
\end{equation}
All together we have found a contour $\mathcal C$ along which the integral (\ref{simpleIntegral}) can be computed with the steepest descent method (see \cite{WongIntegrals}), giving as the following approximation for $\bpar\rightarrow 0$
\begin{equation}
e^{\frac{h(x_0)}{2\pi i \bpar^2N}} \frac{g_{4_1}\left(-\frac23\frac{\pi i}{\sqrt N}\right)}
{\sqrt{iN^{-1}h^\pp (x_0)}}(1+\mathcal O(\bpar^2)
\end{equation}
where 
$$g_{4_1}(x) := \frac{1}{\sqrt N} \sum_{k=0}^{N-1} \phi_{-x}(k)\conj\phi_{x}(k).$$
We now go back to $\chi_{4_1}^{(N)}(0)$, and we  write it as the following integral 
\begin{equation}\label{eq:exacChi41}
\chi_{4_1}^{(N)}(0) =\frac{1}{2\pi \bpar\sqrt{N} }\sum_{k\in\IntN} \int_{\Real - id}f_{\bpar}(x,k)\diff (x,k)
\end{equation}
where 
\begin{equation}
f_\bpar(x,k) = \thDilog\left(\frac{-x}{2\pi \bpar}, - k\right) \conj{\thDilog\left(\frac{x}{2\pi \bpar},-k\right)}.
\end{equation}
Then consider the contour 
\begin{equation}
\mathcal C_{\bpar} = \{z\in\Complex :  \arg f_{\bpar} (z)  = \arg f_{\bpar}(z_{\bpar})\text{, } 
\abs{f_{\bpar} (z)}=\abs{f_{\bpar}(z_{\bpar})}\} 
\end{equation}
where $z_{\bpar}$ is defined as the solution to 
\begin{equation}
\dfield{x}\log f_{\bpar}(x)  = 0
\end{equation}
which minimize the absolute value of $f_{\bpar}$. Using the asymptotic formula for $f_\bpar$ it is simple to show that
the contours $\mathcal C_{\bpar}$ approximates $\mathcal C$ as $\bpar \rightarrow 0$ and that the points $z_{\bpar}$'s will converge to $x_0$. So, in the limit $\bpar \rightarrow 0$, the integral \eqref{eq:exacChi41} is approximated by the integral \eqref{simpleIntegral}, for which we already have an asymptotic formula. We have proved the following
\begin{equation}
\chi_{4_1}^{(N)}(0) =e^{\frac{h(x_0)}{2\pi i \bpar^2N}} \frac{g_{4_1}\left(-\frac23\frac{\pi i}{\sqrt N}\right)}
{\sqrt{iN^{-1}h^\pp (x_0)}}(1+\mathcal O(\bpar^2))\text{,}
\end{equation}
As we remarked above $\im h(x_{0}) = - \operatorname{Vol}(4_1)$.\\
Next we look at the number $g_{4_1}\left(-\frac23\frac{\pi i}{\sqrt N}\right)$ which is a topological invariant of the knot in the formula above. We have that
\begin{align*}
\sqrt Ng_{4_1}&\left(-\frac23\frac{\pi i}{\sqrt N}\right) = \sum_{k=1}^{N} \phi_{\frac23\frac{\pi i}{\sqrt N}}(k)\conj\phi_{-\frac23\frac{\pi i}{\sqrt N}}(k)\\
&=  \abs{\prod_{j= 1}^{N-1} \left( 1 - e^{-\frac{\pi i }{3N}} e^{-\frac{2 \pi i j}{N} }   \right)^{\frac jN }} \sum_{k=0}^{N-1} \prod_{j=1}^{k}\frac{1}{\abs{1-e^{\frac13\frac{-\pi i}{N}}e^{\frac{2\pi i j}{N}} }^2}
\end{align*}
The last expression allows us to make the following remark
\begin{align}\label{eq:41Finite}
g_{4_1}\left(-\frac23\frac{\pi i}{\sqrt N}\right) = \gamma_N \mathcal H^0_N(\conj{\rho_{comp}})
\end{align} 
where $H^0_N(\conj{\rho_{comp}})$ is the Baseilhac--Benedetti invariant for the figure--eight knot found in \cite{BBqhi}, computed at the conjugate of the complete hyperbolic structure (meaning that the holonomies of the structure are all complex conjugated) and $\gamma_N$ is a global rescaling given by
\begin{align}
\gamma_N = \abs{\prod_{j= 1}^{N-1} \left( 1 -  e^{-\frac{2 \pi i j}{N} }   \right)^{\frac jN }}
\end{align}

\begin{remark}
	The very same steps of the previous asymptotic computation for $\chi_{4_1}^{(N)}(0) $ can be applied to $\chi_{5_2}^{(N)}(0)$ up to the point of having an expression
	\begin{align}\label{eq:52asyLog}
	\chi_{5_2}^{(N)}(0) =e^{\frac{\phi(x_{5_2})}{2\pi i \bpar^2N})} \frac{g_{5_1}\left(x_{5_2}\right)}
	{\sqrt{iN^{-1}h_{5_2}^\pp (x_{5_2})}}(1+\mathcal O(\bpar^2))\text{,}
	\end{align}
	where $x_{5_2}$ is the only critical point in the complex plane that contributes to the steepest descent and 
	\begin{align}\label{eq:52Finite}
	g_{5_1}\left(x\right) = \frac{1}{\sqrt N} \sum_{j=0}^{N-1}\conj{\phi_{-x}}(j) \phi_x(j)\conj{\phi_{-x}}(j)\text{.}
	\end{align}
	The fact that $\im \phi(x_{5_2}) = - \operatorname{Vol}(5_2)$, can be seen directly, see for example \cite{AK1}.
	However this situation is already too complicated to allow us to check relations with other theories. The obvious guess is to look for the Baseilhac--Benedetti invariant, but no explicitly computed examples, other then $4_1$, are known to the authors.
\end{remark}

The following conjecture was originally stated in \cite{AK1}  for $N=1$. Here we restate it in the updated setting. 
\begin{conjecture}[\cite{AK1}]\label{cj:VolAK}
	Let $M$ be a closed oriented compact 3-manifold.
	For any hyperbolic knot $K\subset  M$, there exist a two paramters $(\bpar, N)$ family of smooth functions
	$J_{M,K}^{(\bpar,N)}(x,j)$ on $\Real\times\IntN$  which has the following properties.
	\begin{enumerate}
		\item
		For any fully balanced shaped ideal triangulation $X$ of the complement  of $K$ in $M$, there exist a gauge invariant real  linear combination of dihedral angles $\lambda$, a (gauge non-invariant) real quadratic polynomial of dihedral angles $\phi$ such that
		\[
		Z_{\bpar}^{(N)}(X)=e^{i\cb^2{\phi}}\frac{1}{\sqrt N} \sum_{j=0}^{N-1}\int_{\Real} J_{M,K}^{(\bpar,N)}(x, j)e^{i\cb x\lambda} \diff x
		\]
		\item
		For any one vertex shaped H-triangulation $Y$ of the pair $(M,K)$ there exists  a real quadratic polynomial of dihedral angles $\varphi$ such that
		\[
		\lim_{\omega_Y\rightarrow \tau}\bDilog\left(\cb\frac{\omega_Y(K)-\pi}{\pi\sqrt N }, 0\right)Z_{\bpar}^{(N)}(Y)=
		e^{i\cb^2{\varphi}-i\frac{ \pi N}{12}}J_{M,K}^{(\bpar, N)}(0,0),
		\]
		where $\tau\colon \Delta_1(Y)\rightarrow \Real$ takes the value $0$ on the knot $K$ and the value $2\pi$ on all other edges.
		\item The hyperbolic volume of the complement of $K$ in $M$ is recovered as the following limit
		\[
		\lim_{\bpar\rightarrow0}2\pi\bpar^2N\log\vert J_{M,K}^{(\bpar,N)}(0,0)\vert=-\operatorname{Vol}(M\setminus K) 
		\]
	\end{enumerate}
\end{conjecture} 

\begin{remark}\label{rm:AKconj}
	We have proved this extended conjecture for the knots $(S^3,4_1)$ and $(S^3,5_2)$, see formulas  \eqref{eq:41Finite}, \eqref{eq:52asyLog} and \eqref{eq:52Finite}. 
	Moreover we gave a more explicit expansion, showing the appearance of an extra interesting therm $g_K$, and a
	precise relation between $g_{4_1}$ and a known invariant of hyperbolic knots, defined by Baseilhac--Benedetti in \cite{BBqhi}, see equation \eqref{eq:41Finite}.
	We could have been more bold and extend the conjecture declaring the appearance of  $g_{(M,K)}$ to be general, and it to be proportional to the Baseilhac--Benedetti invariant. However we feel that there are not enough evidences to state it as general conjecture.
\end{remark}

\begin{appendices}

	\section{Tempered Distributions}\label{AP:distro}

	For standard references for the topics of this appendix see e.g. \cite{HorBook, HorBook2} and \cite{RSbook2, RSbook1}.

	\begin{definition}
		The Schwartz space $\Schw(\Real^n)$ is the space of all the functions $\phi\in C^{\infty}(\Real^n, \Complex)$ such that
		$$||\phi ||_{\alpha,\beta} :=\sup_{x\in\Real^n} | x^{\beta} \partial^{\alpha} \phi (x)| < \infty$$
		for all multi-indeces $\alpha$, $\beta$.\\
		The space of Tempered Distributions $\Schw'(\Real^n)$ is the space of linear functionals on $\Schw(\Real^n)$ which are continuous with respect to all these seminorms.
	\end{definition}
	Both these spaces are stable under the action of the Fourier transform $\mathcal{F}$ and we use the notation $\hat u = \mathcal{F}(u)$.
	Let $Z_{n}$ be the zero section of $T^*(\Real^n)$.
	\begin{definition}
		For a temperate distribution $u\in \Schw'(\Real^n)$, we define its \emph{Wave Front Set} to be the following subset of the cotangent bundle of $\Real^n$
		$$
		\wf(u) = \{ (x,\xi)\in T^*(\Real^n)- Z_{\Real^n} |\ \xi \in \Sigma_x(u)\}
		$$
		where
		$$
		\Sigma_x(u) = \cap_{\phi\in C^\infty_x(\Real^n)}\Sigma(\phi u).
		$$
		Here
		$$
		C^\infty_x(\Real^n) = \{ \phi \in C^\infty_0(\Real^n) | \phi(x) \neq 0\}
		$$
		and $\Sigma(v)$ are all $\eta \in \Real^n -\{0\}$ having no conic neighborhood $V$ such that
		$$
		|\hat{v}(\xi)| \leq C_N (1+ |\xi|)^{-N}, \ N\in\Integer_{>0},\  \xi\in V.
		$$
		
	\end{definition}

	\begin{lemma}
		Suppose $u$ is a bounded density on a $C^\infty$ sub-manifold $Y$ of $\Real^n$, then $u\in \Schw'(\Real^n)$ and
		$$\wf(u) = \{(x,\xi) \in T^*(\Real^n) | x\in \operatorname{Supp} u\mbox{, } \xi\neq 0 \mbox{ and } \xi(T_xY) = 0 \}.$$
	\end{lemma}
	
	In particular if $\operatorname{Supp} u = Y$, then we see that $\wf(u)$ is the co-normal bundle of $Y$.
	
	\begin{definition}
		Let $u$ and $v$ be temperate distributions on $\Real^n$. Then we define
		$$\wf(u) \oplus \wf(v) = \{ (x,\xi_1+\xi_2) \in T^*(\Real^n) | (x,\xi_1)\in \wf(u)\mbox{, }  (x,\xi_2) \in \wf(v)\}.$$
	\end{definition}
	
	\begin{theorem}
		Let $u$ and $v$ be temperate distributions on $\Real^n$. If
		$$\wf(u) \oplus \wf(v) \cap Z_{n} = \emptyset,$$
		then the product of $u$ and $v$ exists and $uv\in \Schw'(\Real^n)$.
	\end{theorem}
	
	\begin{definition}
		
		We denote by $\Schw(\Real^n)_m$ the set of all $\phi\in C^\infty(\Real^n)$ such that
		$$
		\sup_{x\in \Real^n} | x^\beta \partial^\alpha (\phi)(x)| < \infty
		$$
		for all multi-indices $\alpha$ and $\beta$ such that if $\alpha_i =0$ then $\beta_i=0$ for $n-m<i\leq n$.
		We define $\Schw'(\Real^n)_m$ to be the continuous dual of $\Schw(\Real^n)_m$ with respect to these semi-norms.
	\end{definition}
	
	We observe that if $\pi : \Real^n \longrightarrow \Real^{n-m}$ is the projection onto the first $n-m$ coordinates, then $\pi^*(\Schw(\Real^{n-m})) \subset \Schw(\Real^n)_m$. This means we have a well defined push forward map
	$$
	\pi_* : \Schw'(\Real^n)_m \longrightarrow \Schw'(\Real^{n-m}).
	$$

	\begin{proposition}
		Suppose $Y$ is a linear subspace in $\Real^n$, $u$ a density on $Y$ with exponential decay in all directions in $Y$. Suppose $\pi : \Real^n \longrightarrow \Real^m$ is a projection for some $m<n$. Then $u\in \Schw'(\Real^n)_m$ and $\pi_*(u)$ is a density on $\pi(Y)$ with exponential decay in all directions of the subspace $\pi(Y)\subset \Real^m$.
	\end{proposition}
	
	Tempered distributions can be thought of as functions of growth at most polynomial, thanks to the following 
	\begin{theorem}\label{th:schwPth}
		Let $T\in\Schw'(\Real^n)$, then $T = \partial^{\beta} g$ for some polynomially bounded continuous function $g$ and some multi-index $\beta$. That is, for $f\in\Schw(\Real^n)$,
		$$T (f) = \int_{\Real^n} (-1)^{\abs{\beta}} g(x)(\partial^{\beta} f)(x) \diff x$$
	\end{theorem}
	In particular it is possible to show that $\Schw(\Real^n) \subset \Schw'(\Real^n)$, where $\Schw(\Real^n)\ni f\mapsto T_f\in\Schw'(\Real^n)$ with $T_f(g) = \int_{\Real^n} f(x)g(x)\diff x$.
	
	Denoting by  $ \mathcal{L}( \Schw(\Real^{n}), \Schw'(\Real^{m}))$ the space of continuous linear maps from $\Schw(\Real^{n})$ to $ \Schw'(\Real^{m})$, we remark that we have an isomorphism
	\begin{equation}\label{iso}
	\tilde{\cdot} :   \mathcal{L}( \Schw(\Real^{n}), \Schw'(\Real^{m})) \rightarrow \Schw'(\Real^{n\sqcup m})
	\end{equation}
	determined by the formula
	\begin{equation}\label{eq:NuclearFormula}
	\varphi(f)(g) = \tilde{\varphi}(f\otimes g)
	\end{equation}
	for all $\varphi \in  \mathcal{L}( \Schw(\Real^{n}), \Schw'(\Real^{m}))$, $f\in  \Schw(\Real^{n})$, and $g\in \Schw(\Real^{m})$. This is the content of the Nuclear theorem, see e.g. \cite{RSbook1}. 
	Since we can not freely multiply distributions we  end up with a categoroid instead of a category.
	The partially defined composition in this categroid is  defined as follows. Let $n,m,l$ be three finite sets and $A\in \Schw^\p(\Real^{n\sqcup m})$ and $B\in \Schw^\p(\Real^{m\sqcup l})$. We have pull back maps
	$$
	\pi_{n,m}^* : \Schw'(\Real^{n\sqcup m}) \rightarrow \Schw'(\Real^{n\sqcup m \sqcup l}) \mbox{ and } \pi_{m,l}^* : \Schw'(\Real^{m\sqcup l}) \rightarrow \Schw'(\Real^{n\sqcup m \sqcup l}).
	$$
	By what we summarised above, the product
	$$
	\pi_{n,m}^*(A)\pi_{m,l}^*(B) \in \Schw'(\Real^{n\sqcup m \sqcup l})
	$$
	is well defined provided the wave front sets of $\pi_{n,m}^*(A)$ and $\pi_{m,l}^*(B)$ satisfy the condition
	\begin{equation}\label{wftrans}
	(\wf(\pi_{n,m}^*(A)) \oplus \wf(\pi_{m,l}^*(B)) )\cap Z_{n\sqcup m \sqcup l} = \varnothing
	\end{equation}
	If we now further assume that $\pi_{n,m}^*(A)\pi_{m,l}^*(B)$ continuously extends to $\Schw(\Real^{n\sqcup m \sqcup l})_m$,
	then we obtain a well defined element
	$$
	(\pi_{n,l})_*(\pi_{n,m}^*(A)\pi_{m,l}^*(B) ) \in \Schw'(\Real^{n\sqcup l}).
	$$
	
	\subsection{Bra-Ket Notation}\label{ap:braket}
	We often use the Bra-Ket notation  to make computations with distributions. 
	For $\varphi\in\Schw'(\Real^n)$ a density and $x\in\Real^n$ we will write
	$$\bkt{x|\varphi} :=\varphi(x),$$
	with distributional meaning
	\begin{align*}
	\varphi(f) = \int_{\Real^n}\bkt{x|\varphi} {f(x)}\diff x = \int_{\Real^n}\varphi(x) {f(x)}\diff x\text{.}
	\end{align*}
	In particular if $\varphi\in\Schw(\Real^n)\subset\Schw^\p(\Real^n)$, then
	\begin{align*}
	\bkt{x|\varphi} =\varphi(x) = \delta_x(\varphi)
	\end{align*}
	The integral kernel of the operator $\opT$, if it exists, is a distribution $k_{\opT}$ such that
	\begin{equation}\label{eq:kern1}
	\opT(\psi)(x) = \int_{\Real^{n}}k_{\opT}(x,y)\psi(y)\diff y
	\end{equation}
	Working with Schwartz functions, the nuclear theorem expressed by formula \eqref{eq:NuclearFormula} guarantees that the kernel $k_T$ exists and that it is a tempered distribution. We will usually write the kernel from equation \eqref{eq:kern1}, in Bra-Ket notation as follows
	\begin{equation}
	\opT(\psi)(x) = \int_{\Real^{n}}\bra{x}\opT\ket{y}\psi(y)\diff y
	\end{equation}
	and the nuclear theorem morphism \eqref{eq:NuclearFormula} can be read as
	\begin{equation}\label{eq:NuclearBraKet}
	\bra{x}\opT\ket{y} =\bkt{x,\, y|\widetilde{T}}\text{.}
	\end{equation}
	
	\subsection{$\SqInt(\Affine_N)$ and $\Schw(\Affine_N)$}\label{ap:AffineN}
	$\Affine_N \equiv\Real\times (\IntN)$ has the structure of a locally compact abelian group, with the normalized Haar measure $\diff (x,n)$ defined by
	$$\int_{\Affine_N} f(x,n)\diff (x,n) := \frac{1}{\sqrt{N}}  \sum_{n\in \IntN}\int_{\Real} f(x,n)\diff x$$
	where $f:\Affine_N\longrightarrow \Complex $ is an integrable function.
	By definition $\SqInt(\Affine_N)$ is the space of functions $\mathsf f\colon \Affine_N \longrightarrow \Complex$ such that
	\begin{equation}
	\int_{\Affine_N}\abs{\mathsf{f}(a)}^2\diff a \equiv \frac{1}{\sqrt N}\sum_{n=0}^{N-1} \int_{\Real} \abs{\mathsf f{(x,n)}}^2\diff x <\infty
	\end{equation}
	with standard inner product
	\begin{align}\label{eq:ANinner}
	\bkt{\mathsf f, \mathsf g} \equiv \frac{1}{\sqrt N}\sum_{n=0}^{N-1}\int_{\Real}\mathsf f(x,n)\conj{\mathsf g(x,n)}\diff (x,n)
	\end{align}
	Finite square integrable sequences are just a finite dimensional vector space $$\SqInt(\IntN) \simeq \Complex^N,$$
	with a preferred basis given by mod $N$ Kronecker delta functions
	\begin{align}
	\delta_j(n)  \equiv \left\{ 
	\begin{array}{l l}
	1 & \quad \text{if $j = n$ mod $N$ }\\
	0 & \quad \text{otherwise }\\  
	\end{array} \right. 
	\end{align}
	There is a natural isomorphism
	\begin{align}
	\SqInt(\Real) \otimes \SqInt(\IntN)\simeq\SqInt(\Affine_N)
	\end{align}
	defined by 
	\begin{align}
	f\otimes\delta_j\,(a) = f(x)\delta_j(n)\text{,   for }a=(x,n)\in\Affine_N
	\end{align}
	with inverse
	\begin{align}
	\Affine_N\ni\,\mathsf f\mapsto \sum_{j=0}^{N-1} \mathsf{f} (\cdot, j)\otimes \delta_j\in \SqInt(\Real) \otimes \SqInt(\IntN)
	\end{align}
	Everything just said holds true substituting $\SqInt$ with $\Schw$, with the isomorphism $\Schw(\Affine_N)\simeq\Schw(\Real)\otimes\Complex^N$ and further also, the space of tempered distributions on $\Affine_N$, defined as linear continuous functionals over $\Schw(\Affine_N)$, are simply $\Schw^\p(\Real)\otimes\Complex^N$.
	All the Bra-Ket notation extends trivially to $\Schw(\Affine_N)$, including the nuclear theorem \eqref{eq:NuclearBraKet},  substituting all the integrals over $\Real$ with integrals over $\Affine_N$.\\
	We  use a bracket notation for Fourier coefficients and Gaussian exponentials in $\Affine_N$, following the notation introduced in \cite{AK3}
	\begin{align}\label{eq:FourKern}
	&\bkt{(x,n)\text{, }(y,m)}\equiv e^{2\pi i xy}e^{-2\pi i nm/N}
	&\bkt{(x,n)} \equiv e^{\pi ix^2} e^{-\pi i n(n+N)/N}
	\end{align}
	For $(x,n)$ and $(y,m)$ in $\Affine_N$.
	The Fourier transform then takes the form
	\begin{align*}
	&\mathcal F(f)(x,n) =  \int_{\Affine_N}f(y,m)\bkt{(x,n),(y,m)}\diff (y,m).
	\end{align*}
	For any operator $A$ of order $N$, we can define the operator $\LN (A)$ via the spectral theorem,  such that it formally satisfies 
	$$A = e^{2\pi i \LN (A)/N}\text{.}$$
	We can define, for any function $f:\Affine_N \longrightarrow \Complex$ the operator function $\slashed{f}(\opx,A) \equiv f(\opx,\LN(A))$ for any commuting pair of operators $\opx$ and $A$, where the former is self adjoint and the latter is of order $N$. We have, for $\opx$ and $A$ as above, that
	\begin{equation}\label{eq:FourierSpectra}
	\slashed f(\opx, A) = \int_{\Affine_N}\tilde f(y,m)e^{2\pi iy \opx} A^{-m} \diff (y,m)
	\end{equation}
	where we use the following notation for the inverse Fourier transforms
	\begin{equation}\label{eq:tildeFourier}
	\tilde f (x,n) =\int_{\Affine_N} f(y,m) \conj{\bkt{(y,m);(x,n)}}  \diff (y,m).
	\end{equation}

	\section{Categroids}
	\label{ap:categoroid}
	
	We need a notion which is slightly more general than categories to define the Teichm\"uller TQFT functor.
	
	\begin{definition}\cite{AK1}\\
		A \emph{Categroid} $\mathcal C$ consist of a family of objects $\OBJ({\mathcal C})$ and for any pair of objects $A,B$ from $ \OBJ({\mathcal C})$ a set $\morph_{\mathcal C}(A,B)$ such that
		the following holds
		\begin{description}
			\item[A] For any three objects $A,B,C$ there is a subset $K^{\mathcal C}_{A,B,C} \subset \morph_{\mathcal C}(A,B)\times \morph_{\mathcal C}(B,C)$, called the composable morphisms and a \emph{composition} map
			$$
			\circ : K^{\mathcal C}_{A,B,C} \rightarrow  \morph_{\mathcal C}(A,C).
			$$
			such that composition of composable morphisms is associative.
			\item[B] For any object $A$ we have an identity morphism $1_A \in  \morph_{\mathcal C}(A,A)$ which is composable with any morphism $f\in \morph_{\mathcal C}(A,B)$ or $g\in \morph_{\mathcal C}(B,A)$ and we have the equations
			$$
			1_A\circ f = f \mbox{, and } g \circ 1_A = g.
			$$
		\end{description}
	\end{definition}

\end{appendices}

\noindent 
J{\o}rgen Ellegaard Andersen  and Simone Marzioni, \\
Center for Quantum Geometry of Moduli Spaces,\\
Department of Mathematics,\\ University of Aarhus,\\
DK-8000, Denmark

\label{lastpage}

\begin{thebibliography}{11}

\bibitem[A1]{A1} J.~E. Andersen. 
\newblock  Asymptotic faithfulness of the quantum $SU(n)$
representations of the mapping class groups.  
\newblock {\em Annals of Mathematics.} {\bf
	163}:347--368, 2006.

\bibitem[A2]{A2} J.~E. Andersen.
\newblock The Nielsen-Thurston classification of mapping
classes is determined by TQFT. 
 \newblock{\em J. Math. Kyoto Univ.} {\bf 48}(2):323--338,  2008.

\bibitem[A3]{A4} J.~E. Andersen.  
\newblock The Witten-Reshetikhin-Turaev invariants of finite order mapping tori I. 
\newblock {\em Journal f\"{u}r Reine und Angewandte Mathematik.}  {\bf 681}:1--38, 2013. 

\bibitem[A4]{Aqt}
J.~E. Andersen. 
\newblock Hitchin's connection, {T}oeplitz operators, and symmetry invariant
deformation quantization.
\newblock {\em Quantum Topol.} {\bf 3}(3-4):293--325, 2012.


\bibitem[AG]{AG}
J.~E. Andersen and N.~L. Gammelgaard.
\newblock {The Hitchin-Witten Connection and Complex Quantum Chern-Simons
	Theory}.
\newblock {\em arXiv:1409.1035}, 2014.

\bibitem[AGP]{AGP} J.~E. Andersen, S. Gukov, D. Pei.
\newblock The Verlinde formula for Higgs bundles, 
\newblock {\em arXiv:1608.01761}, 2016.

\bibitem[AH]{AH} J.~E. Andersen \& B. Himpel. 
\newblock The Witten-Reshetikhin-Turaev invariants of finite order mapping tori II
\newblock {\em Quantum Topology.} {\bf 3}:377--421, 2012.

\bibitem[AHJMMc]{AHJMMc} J.~E. Andersen, B.~Himpel, S.~F. J{\o}rgensen, J.~Martens, and B.~McLellan.  
\newblock The {W}itten-{R}eshetikhin-{T}uraev invariant for links in finite
order mapping tori I. 
\newblock {\em Advances in Mathematics.} {\bf 304}:131--178, 2017.

\bibitem[AK1]{AK1}
J.~E. Andersen and R. Kashaev. 
\newblock A {TQFT} from {Q}uantum {T}eichm{\"u}ller {T}heory.
\newblock {\em Comm. Math. Phys.} {\bf 330}(3):887--934, 2014.

\bibitem[AK2]{AK2}
J.~E. Andersen and R. Kashaev.
\newblock A new formulation of the Teichm{\"u}ller TQFT.
\newblock {\em arXiv:1305.4291}, 2013.

\bibitem[AK3]{AK3}
J.~E. Andersen and R. Kashaev.
\newblock {Complex Quantum Chern-Simons}.
\newblock {\em arXiv:1409.1208}, 2014.

\bibitem[AN]{AN} J. E. Andersen and J.-J. K. Nissen, Asymptotic aspects of the Teichm\"{u}ller TQFT, Preprint 2016. 

\bibitem[AU1]{AU1} J.~E. Andersen \& K. Ueno. 
\newblock Abelian Conformal Field theories and
Determinant Bundles.
\newblock {\em International Journal of Mathematics.}  {\bf 18}:919--993, 2007.

\bibitem[AU2]{AU2} J.~E. Andersen \& K. Ueno, 
\newblock Constructing modular functors
from conformal field theories.  
\newblock{\em  Journal of Knot theory and its
Ramifications.} {\bf 16}(2):127--202, 2007.

\bibitem[AU3]{AU3} J.~E. Andersen \& K. Ueno.  
\newblock Modular functors are determined by
their genus zero data.
\newblock{\em  Quantum Topology.} {\bf 3}:255--291, 2012.

\bibitem[AU4]{AU4} J.~E. Andersen \& K. Ueno.  
\newblock Construction of the Witten-Reshetikhin-Turaev TQFT from
conformal field theory.
\newblock{\em Invent. Math.} {\bf 201}(2):519--559, 2015.

\bibitem[ADW]{ADW} S.~Axelrod, S.~Della~Pietra, E.~Witten.  
\newblock Geometric quantization
of Chern Simons gauge theory.
\newblock {\em J.Diff.Geom.} {\bf 33}:787--902, 1991.


\bibitem[BB]{BBqhi}
S. Baseilhac and R. Benedetti. 
\newblock Quantum hyperbolic geometry.
\newblock {\em Algebr. Geom. Topol.} {\bf 7}:845--917, 2007.

\bibitem[B]{B1}C. Blanchet. 
\newblock Hecke algebras, modular categories and
$3$-manifolds quantum invariants.
\newblock{\em Topology.} {\bf 39}(1):193--223, 2000.

\bibitem[BHMV1]{BHMV1} C. Blanchet, N. Habegger, G. Masbaum \&
P. Vogel.  
\newblock Three-manifold invariants derived from the Kauffman Bracket.
\newblock {\em Topology.} {\bf 31}:685--699, 1992.


\bibitem[BHMV2]{BHMV2} C. Blanchet, N. Habegger, G. Masbaum \&
P. Vogel.  
\newblock Topological Quantum Field Theories derived from the
Kauffman bracket. 
\newblock {\em Topology.} {\bf 34}:883--927, 1995.

\bibitem[BMS]{BMS} M. Bordeman, E. Meinrenken \& M.
Schlichenmaier. 
\newblock Toeplitz quantization of K{\"a}hler manifolds and
$gl(N), N \rightarrow \infty$ limit
\newblock {\em Comm. Math. Phys.} {\bf 165}:281--296, 1994.

\bibitem[D]{D3D} T. Dimofte, 
\newblock Complex Chern-Simons theory at level k via the 3d-3d correspondence. 
\newblock {\em Comm. Math. Phys.} {\bf 339}(2):619--662, 2015.

\bibitem[F]{Fqdilog}
L.~D. Faddeev. 
\newblock Discrete {H}eisenberg-{W}eyl group and modular group.
\newblock {\em Lett. Math. Phys.} {\bf 34}(3):249--254, 1995.

\bibitem[FK]{FKqdilog}
L.~D. Faddeev and R.~M. Kashaev.  
\newblock Quantum dilogarithm.
\newblock {\em Modern Phys. Lett. A.} {\bf 9}(5):427--434, 1994.

\bibitem[FKV]{FKV1}
L.~D. Faddeev, R.~M. Kashaev, and A.~Yu. Volkov. 
\newblock Strongly coupled quantum discrete {L}iouville theory. {I}.
{A}lgebraic approach and duality.
\newblock {\em Comm. Math. Phys.} {\bf 219}(1):199--219, 2001.

\bibitem[FG]{FGcluster}
V. Fock and A. Goncharov. 
\newblock Moduli spaces of local systems and higher {T}eichm\"uller theory.
\newblock {\em Publ. Math. Inst. Hautes \'Etudes Sci.} {\bf 103}:1--211, 2006.

\bibitem[FK]{FKmcg}
L. Funar and R.~M. Kashaev. 
\newblock Centrally extended mapping class groups from quantum {T}eichm\"uller
theory.
\newblock {\em Adv. Math.} {\bf 252}:260--291, 2014.


\bibitem[GKT]{GKT}
N. Geer, R. Kashaev, and V. Turaev. 
\newblock Tetrahedral forms in monoidal categories and 3-manifold invariants.
\newblock {\em J. Reine Angew. Math.} {\bf 673}:69--123, 2012.

\bibitem[Hik1]{Hik1}
K. Hikami. 
\newblock Hyperbolicity of partition function and quantum gravity.
\newblock {\em Nuclear Phys. B.} {\bf 616}(3):537--548, 2001.

\bibitem[Hik2]{Hik2}
K. Hikami. 
\newblock Generalized volume conjecture and the {$A$}-polynomials: the
{N}eumann-{Z}agier potential function as a classical limit of the partition
function.
\newblock {\em J. Geom. Phys.} {\bf 57}(9):1895--1940, 2007.

\bibitem[Hit1]{H2} N.~J. Hitchin.  
\newblock The self-duality equations on a Riemann surface. 
\newblock {\em Proc. London Math. Soc.} {\bf  55}(1):59--126, 1987.

\bibitem[Hit2]{H1} N.~J. Hitchin. 
\newblock Flat connections and geometric quantization.
\newblock {\em Comm. Math. Phys.} {\bf 131}:347--380, 1990.

\bibitem[HRS]{HRS}
C.~D. Hodgson, J.~H. Rubinstein, and H. Segerman. 
\newblock Triangulations of hyperbolic 3-manifolds admitting strict angle
structures.
\newblock {\em J. Topol.} {\bf 5}(4):887--908, 2012.

\bibitem[H{\"o}r1]{HorBook2}
L. H{\"o}rmander.  
\newblock Linear partial differential operators.
\newblock {\em Third revised printing. Die Grundlehren der mathematischen
Wissenschaften,} Band {\bf 116}. Springer-Verlag New York Inc., New York, 1969.

\bibitem[H{\"o}r2]{HorBook}
L. H{\"o}rmander. 
\newblock The analysis of linear partial differential operators. {I},
{\em Grundlehren der Mathematischen Wissenschaften [Fundamental
	Principles of Mathematical Sciences]}, {\bf 256}.
\newblock Springer-Verlag, Berlin, second edition, 1990.

\bibitem[KS]{KS}  A. V. Karabegov \& M.  Schlichenmaier. 
\newblock Identification of Berezin-Toeplitz deformation quantization
\newblock {\em J. Reine Angew. Math.} {\bf 540}:49--76,  2001.

\bibitem[K1]{K6j}
R.~M. Kashaev.
\newblock Quantum dilogarithm as a {$6j$}-symbol.
\newblock {\em Modern Phys. Lett. A.} {\bf 9}(40):3757--3768,  1994.

\bibitem[K2]{KVolKnot}
R.~M. Kashaev. 
\newblock The hyperbolic volume of knots from the quantum dilogarithm.
\newblock {\em Lett. Math. Phys.} {\bf 39}(3):269--275, 1997.

\bibitem[K3]{KQuantTeich}
R.~M. Kashaev. 
\newblock Quantization of {T}eichm\"uller spaces and the quantum dilogarithm.
\newblock {\em Lett. Math. Phys.} {\bf 43}(2):105--115, 1998.

\bibitem[K4]{KDehnDil1}
R.~Kashaev.
\newblock The quantum dilogarithm and {D}ehn twists in quantum {T}eichm\"uller
theory.
\newblock In {\em Integrable structures of exactly solvable two-dimensional
	models of quantum field theory.} 
\newblock	 {\em NATO Sci. Ser. II Math. Phys. Chem.} {\bf 35}:211--221. 
Kluwer Acad. Publ., Dordrecht, 2001.


\bibitem[K5]{KPSL2R}
R.~M. Kashaev. 
\newblock Coordinates for the moduli space of flat PSL($2$, $\mathbb{R}$) -
connections.
\newblock {\em Math. Research Letters.} {\bf 12}:23--36, 2005.

\bibitem[K6]{KLiuTeich2}
R.~M. Kashaev. 
\newblock Discrete {L}iouville equation and {T}eichm\"uller theory.
\newblock In {\em Handbook of {T}eichm\"uller theory. {V}olume {III}},
\newblock {\em IRMA Lect. Math. Theor. Phys.} {\bf 17}:821--851. 
Eur. Math.
Soc., Z\"urich, 2012.

\bibitem[Kir]{KirillovDilogIdentities}
A.~N. Kirillov. 
\newblock Dilogarithm identities.
\newblock {\em Progress of Theoretical Physics Supplement.} {\bf 118}:61--142, 1995.

\bibitem[Lac]{Lack}
M. Lackenby. 
\newblock Word hyperbolic {D}ehn surgery.
\newblock {\em Invent. Math.} {\bf 140}(2):243--282, 2000.

\bibitem[Las]{La1} Y. Laszlo.  
\newblock Hitchin's and WZW connections are the
same.
\newblock {\em J. Diff. Geom.} {\bf 49}(3):547--576, 1998.

\bibitem[LT]{LT}
F. Luo and S. Tillmann.
\newblock Angle structures and normal surfaces.
\newblock {\em Trans. Amer. Math. Soc.} {\bf 360}(6):2849--2866,  2008.


\bibitem[M]{M} S. Marzioni.
\newblock Complex Chern-Simons Theory: Knot Invariants and Mapping Class Group Representations.
\newblock {\em PhD Thesis}, Aarhus University, 2016.

\bibitem[MM]{MM}
H. Murakami and J. Murakami. 
\newblock The colored {J}ones polynomials and the simplicial volume of a knot.
\newblock {\em Acta Math.}, {\bf 186}(1):85--104, 2001.

\bibitem[NZ]{NZ}
W.~D. Neumann and Don Zagier.
\newblock Volumes of hyperbolic three-manifolds.
\newblock {\em Topology}, {\bf 24}(3):307--332, 1985.


\bibitem[P]{Pbook}
R.~C. Penner.
\newblock { Decorated {T}eichm\"uller theory}.
\newblock {\em QGM Master Class Series.} European Mathematical Society (EMS),
Z\"urich.
\newblock With a foreword by Yuri I. Manin, 2012.


\bibitem[RS1]{RSbook2}
M. Reed and B. Simon. 
\newblock{Methods of modern mathematical physics. {II}. {F}ourier
	analysis, self-adjointness}.
\newblock Academic Press [Harcourt Brace Jovanovich, Publishers], New
York-London, 1975.


\bibitem[RS2]{RSbook1}
M. Reed and B. Simon. 
\newblock{Methods of modern mathematical physics. {I}}.
\newblock Academic Press, Inc. [Harcourt Brace Jovanovich, Publishers], New
York, second edition, 1980.

\bibitem[RT1]{RT1} N. Reshetikhin \& V. Turaev.
\newblock Ribbon graphs and
their invariants derived from quantum groups
\newblock {\em Comm. Math. Phys.}
{\bf 127}:1--26, 1990.

\bibitem[RT2]{RT2} N. Reshetikhin \& V. Turaev.
\newblock Invariants of
$3$-manifolds via link polynomials and quantum groups
 \newblock {\em Invent. Math.} {\bf 103}:547--597,  1991.


\bibitem[TUY]{TUY} A. Tsuchiya, K. Ueno \& Y. Yamada.
\newblock Conformal Field Theory on
Universal Family of Stable Curves with Gauge Symmetries
\newblock {\em Advanced Studies in	Pure Mathmatics.} {\bf 19}:459--566,  1989.


\bibitem[V]{Vhypergeometry}
A.~Y. Volkov.
\newblock Noncommutative hypergeometry.
\newblock {\em Comm. Math. Phys.} {\bf 258}(2):257--273, 2005.

\bibitem[W1]{W1} E. Witten.
\newblock Quantum field theory and the Jones polynomial.
\newblock {\em Comm. Math. Phys.} {\bf 121}:351--98,  1989.


\bibitem[W2]{W3}
E. Witten.
\newblock Quantization of {C}hern-{S}imons gauge theory with complex gauge
group.
\newblock {\em Comm. Math. Phys.} {\bf 137}(1):29--66, 1991.


\bibitem[Won]{WongIntegrals}
R.~Wong.
\newblock{Asymptotic approximations of integrals}. 
\newblock {\em
	Classics in Applied Mathematics}, {\bf 34}.
\newblock Society for Industrial and Applied Mathematics (SIAM), Philadelphia,
PA.
\newblock Corrected reprint of the 1989 original, 2001.

\end{thebibliography}
\end{document}